\documentclass[a4paper,leqno]{article}

\usepackage{ifthen}
\usepackage[latin1]{inputenc}
\usepackage[T1]{fontenc}
\usepackage{xspace}
\usepackage{amsfonts,amssymb,amsxtra}
\usepackage{amsmath}
\usepackage{amsthm}
\usepackage[all]{xy}
\usepackage[english,french]{babel}

\theoremstyle{definition}
\newtheorem{ssect}[subsection]{}
\newtheorem{defn}[subsection]{Définition}
\newtheorem*{remqe}{Remarque}
\newtheorem{remq}[subsection]{Remarque}

\newtheorem{vari}[subsection]{Variante}
\theoremstyle{plain}
\newtheorem{lemm}[subsection]{Lemme}
\newtheorem{prop}[subsection]{Proposition}
\newtheorem{theo}[subsection]{Théorème}
\newtheorem{coro}[subsection]{Corollaire}

\newcommand{\beq}{\begin{equation}}
\newcommand{\eeq}{\end{equation}}

\newcommand{\al}{\alpha}

\newcommand{\eps}{\epsilon}
\newcommand{\Lbd}{\Lambda}
\newcommand{\lbd}{\lambda}

\newcommand{\cF}{\mathcal{F}}
\newcommand{\cG}{\mathcal{G}}
\newcommand{\fH}{\mathcal{H}}
\newcommand{\cO}{\mathcal{O}}

\newcommand{\fm}{\mathfrak{m}}
\newcommand{\bA}{\mathbb{A}}
\newcommand{\bC}{\mathbb{C}}
\newcommand{\F}{\mathbb{F}}
\newcommand{\N}{\mathbb{N}}
\newcommand{\Q}{\mathbb{Q}}
\newcommand{\R}{\mathbb{R}}
\newcommand{\Z}{\mathbb{Z}}

\newcommand{\lbar}{\overline}

\newcommand{\ot}{\otimes}
\newcommand{\hra}{\hookrightarrow}

\newcommand{\op}{^\circ}

\newcommand{\Qb}{\overline{\Q}}
\newcommand{\Qlb}{\overline{\Q_l}}
\newcommand{\Dbc}{D^b_c}
\newcommand{\idI}{{i \in I}}

\newcommand{\Fb}{\overline{F}}

\newcommand{\Kb}{\overline{K}}

\newcommand{\xb}{{\bar{x}}}
\newcommand{\yb}{{\bar{y}}}
\newcommand{\kb}{\bar{k}}
\newcommand{\sbar}{{\bar{s}}}
\newcommand{\Sbar}{{\overline{S}}}
\newcommand{\etab}{{\bar{\eta}}}
\newcommand{\Zch}{\hat{\Z}'}
\newcommand{\simto}[1][]{\sto[#1]{\sim}}
\newcommand{\pimod}{\pi_1^{\mathrm{mod}}}
\newcommand{\ii}{\iota_i}

\newcommand{\sm}{_\bullet}
\newcommand{\ip}{{(i)}}
\newcommand{\jp}{{(j)}}
\newcommand{\ijp}{{(ij)}}
\newcommand{\np}{{(n)}}
\newcommand{\xp}{{(x)}}

\newcommand{\cGc}{\check{\cG}}
\newcommand{\hookto}[1]{\stackrel{#1}{\hookrightarrow}}

\newcommand{\DiUi}{{D_{i,U_\ip}}}
\newcommand{\DjUj}{{D_{j,U_\jp}}}
\newcommand{\cHom}{\fH\mathit{om}}
\newcommand{\Xt}{\widetilde{X}}

\newcommand{\ensdr}[2]{\left\{ #1 \,\left|\, #2 \right.\right\}}

\newcommand{\sto}[2][]{
\xrightarrow[#1]{#2}} 

\newcommand{\Hom}{\mathrm{Hom}}

\newcommand{\GL}{\mathrm{GL}}
\newcommand{\Mod}{\mathrm{Mod}}

\newcommand{\Reg}{\mathrm{Reg}}
\newcommand{\gr}{\mathrm{gr}}
\newcommand{\cosq}{\mathrm{cosq}}
\newcommand{\id}{\mathrm{id}}

\newcommand{\car}{\mathrm{car}}
\newcommand{\red}{\mathrm{red}}
\newcommand{\codim}{\mathrm{codim}}

\newcommand{\Modcent}[2][]{\Mod_c(#2)_{#1\textrm{\ifthenelse{\equal{#1}{}}{}{-}}\mathrm{ent}}}
\newcommand{\Modcentinv}[2][]{\Mod_c(#2)_{#1\textrm{\ifthenelse{\equal{#1}{}}{}{-}}\mathrm{ent}^{-1}}}
\newcommand{\Modclent}[2][]{\Modcent[#1]{#2,\Qlb}}
\newcommand{\Modclentinv}[2][]{\Modcentinv[#1]{#2,\Qlb}}
\newcommand{\Dbcent}[2][]{\Dbc(#2)_{\ifthenelse{\equal{#1}{}}{}{#1\textrm{-}}\mathrm{ent}}}
\newcommand{\Dbceent}[2][]{\Dbc(#2)_{\ifthenelse{\equal{#1}{}}{\eps}{\deeps
#1.}\textrm{-}\mathrm{ent}}}
\newcommand{\Dbcentinv}[2][]{\Dbc(#2)_{\ifthenelse{\equal{#1}{}}{}{#1\textrm{-}}\mathrm{ent}^{-1}}}
\newcommand{\Dbceentinv}[2][]{\Dbc(#2)_{\ifthenelse{\equal{#1}{}}{\eps}{\deeps
#1.}\textrm{-}\mathrm{ent}^{-1}}}
\def\deeps#1#2#3#4.{\ifthenelse{\equal{#3}{0}}{I}{(I\ifthenelse{\equal{#3}{-}}{}{+}#3#4)}}
\newcommand{\Dbclent}[2][]{\Dbcent[#1]{#2,\Qlb}}
\newcommand{\Dbcleent}[2][]{\Dbceent[#1]{#2,\Qlb}}
\newcommand{\Dbclentinv}[2][]{\Dbcentinv[#1]{#2,\Qlb}}
\newcommand{\Dbcleentinv}[2][]{\Dbceentinv[#1]{#2,\Qlb}}

\newcommand{\dModcent}[3][]{\Modcent[#1]{#2}^{#3}}

\newcommand{\dModclent}[3][]{\Modclent[#1]{#2}^{#3}}
\newcommand{\dModclentinv}[3][]{\Modclentinv[#1]{#2}^{#3}}
\newcommand{\dDbcent}[3][]{{\Dbcent[#1]{#2}^{#3}}}
\newcommand{\dDbceent}[3][]{{\Dbceent[#1]{#2}^{#3}}}
\newcommand{\dDbcentinv}[3][]{{\Dbcentinv[#1]{#2}^{#3}}}
\newcommand{\dDbceentinv}[3][]{{\Dbceentinv[#1]{#2}^{#3}}}
\newcommand{\dDbclent}[3][]{\Dbclent[#1]{#2}^{#3}}
\newcommand{\dDbcleent}[3][]{\Dbcleent[#1]{#2}^{#3}}
\newcommand{\dDbclentinv}[3][]{\Dbclentinv[#1]{#2}^{#3}}
\newcommand{\dDbcleentinv}[3][]{\Dbcleentinv[#1]{#2}^{#3}}

\DeclareMathOperator{\Spec}{Spec}
\DeclareMathOperator{\Tr}{Tr}

\DeclareMathOperator{\Gal}{Gal}

\DeclareMathOperator{\Ker}{Ker}

\newcommand{\tm}[1]{\textnormal{(#1)}}
\newcommand{\tmi}{\tm{i} }
\newcommand{\tmii}{\tm{ii} }
\newcommand{\tmiii}{\tm{iii} }

\newcommand{\bieqdim}{biéquidimensionnel\xspace}
\newcommand{\cacart}{carré cartésien\xspace}

\newcommand{\can}{canonique\xspace}
\newcommand{\cara}{caractéristique\xspace}

\newcommand{\cat}{catégorie\xspace}

\newcommand{\clot}{clôture\xspace}

\newcommand{\comms}{commutatifs\xspace}
\newcommand{\comm}{commutatif\xspace}
\newcommand{\compl}{complémentaire\xspace}

\newcommand{\cst}{constant\xspace}

\newcommand{\dem}{démonstration\xspace}
\newcommand{\Demo}{Démonstration\xspace}

\newcommand{\diag}{diagramme\xspace}
\newcommand{\diacom}{diagramme commutatif\xspace}
\newcommand{\dist}{distingué\xspace}
\newcommand{\Dpr}{D'après }
\newcommand{\dpr}{d'après }

\newcommand{\elt}{élément\xspace}
\newcommand{\Enpart}{En particulier\xspace}

\newcommand{\etale}{étale\xspace}

\newcommand{\faisc}{faisceau\xspace}
\newcommand{\faisx}{faisceaux\xspace}

\newcommand{\gel}{général\xspace}

\newcommand{\geom}{géométrique\xspace}
\newcommand{\geot}{géométriquement\xspace}

\newcommand{\hyperr}{hyperrecouvrement\xspace}

\newcommand{\hyp}{hypothèse\xspace}
\newcommand{\ie}{\emph{i.~e.}, }
\newcommand{\inte}{intégralité\xspace}
\newcommand{\irrs}{irréductibles\xspace}
\newcommand{\irr}{irréductible\xspace}
\newcommand{\iso}{isomorphisme\xspace}
\newcommand{\isoms}{isomorphismes\xspace}
\newcommand{\isom}{isomorphisme\xspace}

\newcommand{\loccst}{\loct \cst}
\newcommand{\loct}{localement\xspace}
\newcommand{\me}{même\xspace}

\newcommand{\mrs}{modérément ramifiés\xspace}
\newcommand{\mr}{modérément ramifié\xspace}
\newcommand{\noeths}{noethériens\xspace}
\newcommand{\noeth}{noethérien\xspace}
\newcommand{\Ops}{On peut supposer\xspace}
\newcommand{\ops}{on peut supposer\xspace}

\newcommand{\rec}{récurrence\xspace}
\newcommand{\reg}{régulier\xspace}

\newcommand{\relta}{\relt à\xspace}
\newcommand{\relt}{relativement\xspace}

\newcommand{\repn}{représentation\xspace}

\newcommand{\resp}{resp.\@\xspace} 
\newcommand{\Resp}{Resp.\@\xspace}

\newcommand{\rev}{revêtement\xspace}
\newcommand{\schs}{schémas\xspace}
\newcommand{\sch}{schéma\xspace}

\newcommand{\sepble}{séparable\xspace}

\newcommand{\sep}{séparé\xspace}

\newcommand{\ssi}{si et seulement si\xspace}

\newcommand{\ssch}{sous-schéma\xspace}

\newcommand{\teq}{telle que\xspace}

\newcommand{\tf}{de type fini\xspace}

\newcommand{\thm}{théorème\xspace}
\newcommand{\topet}{\topo \etale}
\newcommand{\topo}{topologie\xspace}

\newcommand{\tq}{tel que\xspace}
\newcommand{\trdist}{triangle distingué\xspace}

\newcommand{\tsq}{tels que\xspace}
\newcommand{\unipe}{unipotente\xspace}

\newcommand{\cn}{à croisements normaux\xspace}
\newcommand{\scn}{strictement \cn}
\newcommand{\dcn}{diviseur \cn}
\newcommand{\dscn}{diviseur \scn}
\newcommand{\sst}{semi-stable\xspace}
\newcommand{\sss}{strictement semi-stable\xspace}

\newcommand{\ext}{extension\xspace}
\newcommand{\seps}{séparés\xspace}
\newcommand{\seult}{seulement\xspace}

\newcommand{\surj}{surjectif\xspace}
\newcommand{\cb}{changement de base\xspace}
\newcommand{\morph}{morphisme\xspace}
\newcommand{\regs}{réguliers\xspace}
\newcommand{\constr}{constructible\xspace}
\newcommand{\csts}{constants\xspace}
\newcommand{\cic}{complètement \ic}
\newcommand{\ic}{intégralement clos\xspace}
\newcommand{\propo}{proposition\xspace}

\newcommand{\chgt}{changement\xspace}

\newcommand{\ver}{vérifi}

\newcommand{\ladiq}{\hbox{$l$-adique}\xspace}

\newcommand{\repnl}{\repn \ladiq}
\newcommand{\ntron}{\hbox{$n$-tronqué}\xspace}
\newcommand{\sscin}{\hbox{$s$-scindé}\xspace}
\newcommand{\sscins}{\hbox{$s$-scindés}\xspace}

\newcommand{\eent}[1][]{$\ifthenelse{\equal{#1}{}}{\eps}{\deeps
#1.}$-entier\xspace}
\newcommand{\eentinv}[1][]{$\ifthenelse{\equal{#1}{}}{\eps}{\deeps
#1.}$-entier inverse\xspace}
\newcommand{\ent}[1][]{\ifthenelse{\equal{#1}{}}{}{$#1$-}entier\xspace}
\newcommand{\entinv}[1][]{\ifthenelse{\equal{#1}{}}{}{$#1$-}entier inverse\xspace}
\newcommand{\eents}[1][]{$\ifthenelse{\equal{#1}{}}{\eps}{\deeps
#1.}$-entiers\xspace}
\newcommand{\eentinvs}[1][]{$\ifthenelse{\equal{#1}{}}{\eps}{\deeps
#1.}$-entiers inverses\xspace}
\newcommand{\ents}[1][]{\ifthenelse{\equal{#1}{}}{}{$#1$-}entiers\xspace}
\newcommand{\entinvs}[1][]{\ifthenelse{\equal{#1}{}}{}{$#1$-}entiers inverses\xspace}

\newcommand{\ibid}[1][]{[\emph{ibid.}\ifthenelse{\equal{#1}{}}{}{, }#1]\xspace} 

\newenvironment{cond}{\begin{list}{}{%
\setlength{\leftmargin}{0mm}%
\setlength{\rightmargin}{0mm}%
\setlength{\itemindent}{\parindent}%
\setlength{\listparindent}{\parindent}%
\setlength{\parsep}{\parskip}%
}\item[]\ignorespaces}{\unskip\end{list}}

\makeatletter
\newenvironment{numcond}[1][subsubsection]{\stepcounter{#1}%
\begin{list}{}{%
\setlength{\leftmargin}{0mm}%
\setlength{\rightmargin}{0mm}%
\setlength{\labelwidth}{0mm}%
\setlength{\labelsep}{\parindent}%
\setlength{\itemindent}{\parindent}%
\setlength{\listparindent}{\parindent}%
\setlength{\parsep}{\parskip}%
}\def\@currentlabel{\csname the#1\endcsname}%
\item[\textbf{\@currentlabel}]\ignorespaces}{\unskip\end{list}}
\makeatother

\newcommand{\excludeversion}[1]{\newboolean{#1}\setboolean{#1}{false}}
\newcommand{\vcb}[2]%
 {\ifthenelse{\boolean{#1}}{\cbstart}{}#2\ifthenelse{\boolean{#1}}{\cbend}{}}

\makeatletter
\newcommand{\sigle}[1]{\protect\@sigle{#1}}
\newcommand{\@sigle}[1]{\gdef\@sigleparam{#1}\aftergroup\@absorbtwo}
\newcommand{\@absorbtwo}[2]{\ifcat a#1#1#2\else \@sigleparam \fi}
\makeatother

\excludeversion{v1m} \excludeversion{v1n} \emergencystretch=1em
\usepackage{aeguill}  \DeclareMathSymbol{A}{\mathalpha}{operators}{"41}
  \DeclareMathSymbol{B}{\mathalpha}{operators}{"42}
  \DeclareMathSymbol{C}{\mathalpha}{operators}{"43}
  \DeclareMathSymbol{D}{\mathalpha}{operators}{"44}
  \DeclareMathSymbol{E}{\mathalpha}{operators}{"45}
  \DeclareMathSymbol{F}{\mathalpha}{operators}{"46}
  \DeclareMathSymbol{G}{\mathalpha}{operators}{"47}
  \DeclareMathSymbol{H}{\mathalpha}{operators}{"48}
  \DeclareMathSymbol{I}{\mathalpha}{operators}{"49}
  \DeclareMathSymbol{J}{\mathalpha}{operators}{"4A}
  \DeclareMathSymbol{K}{\mathalpha}{operators}{"4B}
  \DeclareMathSymbol{L}{\mathalpha}{operators}{"4C}
  \DeclareMathSymbol{M}{\mathalpha}{operators}{"4D}
  \DeclareMathSymbol{N}{\mathalpha}{operators}{"4E}
  \DeclareMathSymbol{O}{\mathalpha}{operators}{"4F}
  \DeclareMathSymbol{P}{\mathalpha}{operators}{"50}
  \DeclareMathSymbol{Q}{\mathalpha}{operators}{"51}
  \DeclareMathSymbol{R}{\mathalpha}{operators}{"52}
  \DeclareMathSymbol{S}{\mathalpha}{operators}{"53}
  \DeclareMathSymbol{T}{\mathalpha}{operators}{"54}
  \DeclareMathSymbol{U}{\mathalpha}{operators}{"55}
  \DeclareMathSymbol{V}{\mathalpha}{operators}{"56}
  \DeclareMathSymbol{W}{\mathalpha}{operators}{"57}
  \DeclareMathSymbol{X}{\mathalpha}{operators}{"58}
  \DeclareMathSymbol{Y}{\mathalpha}{operators}{"59}
  \DeclareMathSymbol{Z}{\mathalpha}{operators}{"5A}

\DeclareSymbolFont{psym}{U}{psy}{m}{n}

\DeclareMathSymbol{\alpha}{\mathord}{psym}{97}
\DeclareMathSymbol{\beta}{\mathord}{psym}{98}
\DeclareMathSymbol{\gamma}{\mathord}{psym}{103}
\DeclareMathSymbol{\delta}{\mathord}{psym}{100}
\DeclareMathSymbol{\epsilon}{\mathord}{psym}{101}
\DeclareMathSymbol{\varepsilon}{\mathord}{psym}{101}
\DeclareMathSymbol{\zeta}{\mathord}{psym}{122}
\DeclareMathSymbol{\eta}{\mathord}{psym}{104}
\DeclareMathSymbol{\theta}{\mathord}{psym}{113}
\DeclareMathSymbol{\vartheta}{\mathord}{psym}{74}
\DeclareMathSymbol{\iota}{\mathord}{psym}{105}
\DeclareMathSymbol{\kappa}{\mathord}{psym}{107}
\DeclareMathSymbol{\lambda}{\mathord}{psym}{108}
\DeclareMathSymbol{\mu}{\mathord}{psym}{109}
\DeclareMathSymbol{\nu}{\mathord}{psym}{110}
\DeclareMathSymbol{\xi}{\mathord}{psym}{120}
\DeclareMathSymbol{\omicron}{\mathord}{psym}{111}
\DeclareMathSymbol{\pi}{\mathord}{psym}{112}
\DeclareMathSymbol{\varpi}{\mathord}{psym}{118}
\DeclareMathSymbol{\rho}{\mathord}{psym}{114}
\DeclareMathSymbol{\sigma}{\mathord}{psym}{115}
\DeclareMathSymbol{\varsigma}{\mathord}{psym}{86}
\DeclareMathSymbol{\tau}{\mathord}{psym}{116}
\DeclareMathSymbol{\upsilon}{\mathord}{psym}{117}
\DeclareMathSymbol{\phi}{\mathord}{psym}{102}
\DeclareMathSymbol{\varphi}{\mathord}{psym}{106}
\DeclareMathSymbol{\chi}{\mathord}{psym}{99}
\DeclareMathSymbol{\psi}{\mathord}{psym}{121}
\DeclareMathSymbol{\omega}{\mathord}{psym}{119}

\usepackage[unicode]{hyperref}
\numberwithin{equation}{subsection}

\begin{document}

\author{Weizhe Zheng}
\title{Sur la cohomologie des faisceaux $l$-adiques entiers sur les corps locaux}
\date{}

\maketitle

{\renewcommand{\thefootnote}{}
\footnotetext{W. Zheng, Université Paris-Sud 11, Mathématiques,
Bât.\ 425, 91405 Orsay Cedex, France. Courriel :
\texttt{weizhe.zheng@math.u-psud.fr}} \footnotetext{Classification
mathématique par sujets (2000) : 14F20, 14G20, 11G25, 14D05.}
\footnotetext{Mots clefs :
intégralité, cohomologie $l$-adique, cycles proches.}%
}

\begin{abstract}
On étudie le comportement des faisceaux $l$-adiques entiers sur les
schémas de type fini sur un corps local par les six opérations et
le foncteur des cycles proches.
\end{abstract}


\section{Introduction}\label{S1}

Soient $R$ un anneau de valuation discrète hensélien excellent de
corps résiduel fini de caractéristique~$p$, $K$ son corps des
fractions. Un tel corps sera appelé \emph{corps local}. Soit $\eta
= \Spec K$.

Soit $X$ un \sch \tf sur $\eta$. On désigne par $|X|$ l'ensemble de
ses points fermés. Pour $x\in |X|$, le corps résiduel $\kappa(x)$
de $X$ en $x$ est une extension finie de~$K$. On note $R_x$ son
anneau des entiers, $x_0$ le point fermé de $\Spec R_x$. Soient
$\xb$ un point \geom de $X$ au-dessus de $x$ de corps résiduel
$\kappa(\xb)$ une clôture séparable de $\kappa(x)$, $R_\xb$ la
normalisation de $R_x$ dans $\kappa(\xb)$, $\xb_0$ le point fermé
de $\Spec R_\xb$. Soit $F_x \in \Gal(\kappa(\xb_0)/\kappa(x_0))$ le
Frobenius \geom qui envoie $a$ sur $a^{1/q}$, où $q=\sharp
\kappa(x_0)$.

Fixons un nombre premier $l \neq p$. On désigne par $\Qlb$ une
clôture algébrique de $\Q_l$. Soit $\cF$ un $\Qlb$-\faisc sur $X$.
D'après le \thm de monodromie locale, les valeurs propres d'un
relèvement $\Phi_x\in \Gal(\kappa(\xb)/\kappa(x))$ de $F_x$ agissant
sur $\cF_\xb$ sont bien définies à multiplication près par des
racines de l'unité \cite[1.7.4]{WeilII}.

Rappelons qu'on dit que $\cF$ est \emph{entier} \cite[0.1]{DE} si
les valeurs propres de $\Phi_x$ sont des entiers algébriques pour
tout $x\in |X|$. Cette intégralité est stable par image directe à
support propre \ibid[0.2]. La démonstration utilise l'analogue de ce
résultat sur un corps fini \cite[XXI 5.2.2]{SGA7}.

L'objet de cet article est d'étudier, plus généralement, le
comportement de l'intégralité par les foncteurs usuels : les six
opérations et le foncteur des cycles proches. De façon plus précise,
on examine le comportement par ces foncteurs de la divisibilité des
valeurs propres des $\Phi_x$ par des puissances de $q$. On introduit
pour cela une mesure de la $q$-divisibilité inspirée des \og jauges
\fg de Mazur-Ogus. On prouve notamment les résultats espérés dans
\cite[5.5]{IllMisc}.

\vcb{v1m}{Dans un travail ultérieur \cite{Zheng}, on examine le
comportement de la rationalité et de l'indépendance de $l$ par les
mêmes opérations.}

Les résultats concernant les six opérations sont exposés au
\S~\ref{S2}. Au \S~\ref{S3} on traite le cas crucial de $Rj_*\cF$,
pour l'inclusion $j: U\to X$ du complémentaire d'un diviseur à
croisements normaux $D$ dans un schéma $X$ lisse sur $\eta$ et d'un
faisceau $\cF$ lisse sur $U$ et modérément ramifié le long
de~$D$. Les démonstrations des résultats du \S~\ref{S2} sont
données au \S~\ref{S4}. L'ingrédient essentiel est un théorème
de de Jong, grâce auquel on se réduit au cas traité au
\S~\ref{S3} par les techniques usuelles de descente cohomologique.
Le résultat principal du \S~\ref{S5} est la stabilité de l'\inte
par le foncteur des cycles proches $R\Psi$. À nouveau,
l'ingrédient clef est un \thm de de Jong, qui permet de se ramener
au cas d'un couple \sss et d'un faisceau lisse sur le
complémentaire du diviseur $D$ réunion de la fibre spéciale et
des composantes horizontales et modérément ramifié le long
de~$D$. L'étude de ce cas, plus délicate qu'on ne pouvait s'y
attendre, repose sur une compatibilité technique (\ref{lemm.RPdev}
(ii)) généralisant \cite[1.5 (a)]{IllPL}. Au \S~\ref{S6} on
généralise la notion d'intégralité aux champs algébriques.

\vcb{v1n}{Je remercie chaleureusement L.~Illusie pour m'avoir
suggéré ce sujet, pour son aide à la composition de cet article,
et pour sa lecture minutieuse des diverses versions du manuscrit. Je
suis reconnaissant à G.~Laumon pour une simplification de la
démonstration de \ref{lemm.RPdev} (ii). Je remercie également
O.~Gabber, F.~Orgogozo et le rapporteur pour leurs remarques et
suggestions.}

\section{Intégralité et six opérations}\label{S2}

On conserve les notations du \S~\ref{S1}. On désigne par $\Qb$ la
clôture algébrique de $\Q$ dans~$\bC$. Pour $r\in \Q$, on note
$q^r$ l'unique élément de $\Qb\cap \R_{>0}$ vérifiant $(q^r)^b =
q^a$, où $a, b \in \Z$ sont tel que $r=\frac{a}{b}$, $b\neq 0$.
Soit $X$ un \sch \tf sur $\eta$.

\begin{defn}\label{defn.f}
Fixons un plongement  $\iota : \Qb \to \Qlb$. On dit qu'un
$\Qlb$-\faisc $\cF$ sur $X$ est $r$-\emph{entier} (\resp
$r$-\emph{\entinv{}}) si pour tout point fermé $x$ de $X$, et toute
valeur propre $\al$ de $\Phi_x$ agissant sur $\cF_\xb$,
$\al/\iota(q^r)$ (\resp $\iota(q^r)/\al$) est entier sur $\Z$, où
$q=\sharp \kappa(x_0)$. Cette définition ne dépend pas des choix
de $\Phi_x$ et de $\iota$. On dit que $\cF$ est \emph{entier} (\resp
\emph{entier inverse}) s'il est $0$-entier (\resp $0$-\entinv).
\end{defn}

Les $\Qlb$-\faisx \ents (\resp $r$-entiers, \resp \entinvs, \resp
$r$-\entinvs) sur~$X$ forment une sous-\cat épaisse
\cite[1.11]{Tohoku} de $\Mod_c(X,\Qlb)$, notée $\Modclent{X}$
(\resp $\Modclent[r]{X}$,
\resp $\Modclentinv{X}$, \resp $\Modclentinv[r]{X}$). 

\vcb{v1m}{Soient $K'$ une \ext finie de $K$, $Z$ un \sch \tf sur
$K'$, $\cG\in \Mod_c(Z,\Qlb)$. Alors $\cG$ est $r$-entier (\resp
$r$-entier inverse) \relta $K'$ si et \seult s'il est $r$-entier
(\resp $r$-entier inverse) \relta $K$.}

Rappelons que pour les schémas $X$ séparés de type fini sur un \sch
$S$ régulier de dimension $\le 1$, et en particulier sur $\eta$, on
dispose, par \cite[\S~6]{Ekedahl}, d'une catégorie triangulée
$\Dbc(X,\Qlb)$ et d'un formalisme de six opérations : $Rf_*$,
$Rf_!$, $f^*$, $Rf^!$, $\ot$, $R\cHom$. La catégorie $\Dbc(X,\Qlb)$
est la $2$-limite inductive des catégories $\Dbc(X,E_\lbd)$, où
$E_\lbd$ parcourt les extensions finies de $\Q_l$ contenues dans
$\Qlb$. Si $\cO_\lbd$ est l'anneau des entiers de $E_\lbd$,
$\Dbc(X,E_\lbd)$ est déduite de la catégorie $\Dbc(X,\cO_\lbd)$
définie dans \ibid par extension des scalaires de $\cO_\lbd$ à
$E_\lbd$. Le formalisme construit dans \ibid pour $\Dbc(-,
\cO_\lbd)$ se transpose trivialement.

Ce formalisme a un sens pour les \schs \tf sur $S$ (pas
nécessairement \seps), et ce n'est que pour certaines opérations
($Rf_!$ et $Rf^!$) qu'on a besoin d'une \hyp de séparation sur les
morphismes. Pour un formalisme sans hypothèse de séparation, voir
l'appendice (\S~\ref{S6}).

La définition qui suit est inspirée de la notion des \og jauges \fg
de Mazur-Ogus \cite[8.7]{BO}.

\begin{defn}\label{defn.c}
Soit $\eps: \Z \to \Q$ une fonction. On dit qu'un objet $K\in
\Dbc(X,\Qlb)$ est \emph{entier} (\resp $\eps$-\emph{entier}, \resp
\emph{entier inverse}, \resp $\eps$-\emph{\entinv{}}) si pour tout
$i\in \Z$, $\fH^i(K)$ est entier (\resp $\eps(i)$-entier, \resp
\entinv, \resp $\eps(i)$-\entinv).
\end{defn}

On désigne la sous-\cat pleine de $\Dbc(X,\Qlb)$ formée des objets
\ents (\resp \eents, \resp \entinvs, \resp \eentinvs) par
\[\text{$\Dbclent{X}$ (\resp $\Dbcleent{X}$, \resp $\Dbclentinv{X}$,
\resp $\Dbcleentinv{X}$).}
\]
Lorsque $\eps$ est constant, $\Dbcleent{X}$ et $\Dbcleentinv{X}$
sont des sous-catégories triangulées. On abrège parfois
$\Dbc(X,\Qlb)$ en $\Dbc(X)$.

On note $I$ la fonction d'inclusion de $\Z$ dans $\Q$.

\begin{ssect}
Soient $r$, $r_1$, $r_2\in \Q$.

Pour $\cF \in\Modcent[r_1]{X,\Qlb}$, $\cG\in
\Modcent[r_2]{X,\Qlb}$, on a \[\cF\otimes \cG \in
\Modcent[(r_1+r_2)]{X,\Qlb}.\] Pour $K\in
\Dbcent[(rI+r_1)]{X,\Qlb}$, $L\in \Dbcent[(rI+r_2)]{X,\Qlb}$, on a
\[K\otimes L \in \Dbcent[(rI+r_1+r_2)]{X,\Qlb}.\] De \me pour \og
entier inverse \fg.

Pour $\cF \in\Modclentinv[r_1]{X}$ lisse, $\cG\in
\Modclent[r_2]{X}$, on a \[\cHom(\cF,\cG) \in
\Modclent[(r_2-r_1)]{X}.\] Pour $\cF \in\Modclent[r_1]{X}$ lisse,
$\cG\in \Modclentinv[r_2]{X}$, on a
\[\cHom(\cF,\cG) \in \Modclentinv[(r_2-r_1)]{X}.\]

Soit $f: X\to Y$ un morphisme de \schs \tf sur $\eta$. Alors $f^*$
préserve les complexes $\eps$-entiers (\resp $\eps$-entiers
inverses).
\end{ssect}

\begin{theo}\label{theo.Rf!b} Soient $f: X\to Y$ un morphisme \sep de \schs \tf
sur~$\eta$, $\cF$  un $\Qlb$-\faisc entier (\resp \entinv) sur $X$.
Alors pour tout point fermé $y$ de~$Y$, $(Rf_! \cF)_y$ est entier et
$(I-n)$-entier (\resp $I$-\entinv et $n$-\entinv), où
$n=\dim(f^{-1}(y))$. \Enpart, $Rf_!$ induit
\begin{align}
\Dbcent{X} &\to \Dbcent{Y},\label{eq.Rf!b1}
\\
\Dbceent[1,0]{X} &\to \Dbceent[1,-d_r]{Y},\label{eq.Rf!b2}
\\
\Dbceentinv[1,0]{X} &\to \Dbceentinv[1,0]{Y},\label{eq.Rf!b3}
\\
\Dbcentinv{X} &\to \Dbcentinv[d_r]{Y},\label{eq.Rf!b4}
\end{align}
où $d_r = \max_{y\in |Y|} \dim f^{-1}(y)$ est la dimension
relative.
\end{theo}

Le cas \og entier \fg (\eqref{eq.Rf!b1} et \eqref{eq.Rf!b2}) de
\ref{theo.Rf!b} est un théorème de Deligne-Esnault \cite[0.2]{DE}.

\begin{theo}\label{theo.Rf*b}
Soient $f: X\to Y$ un morphisme \sep de \schs \tf sur $\eta$,
$d_X=\dim X$. Alors $Rf_*$ induit
\begin{align}
\Dbcent{X} &\to \Dbcent{Y},\label{eq.Rf*1}\\
\Dbceent[1,0]{X} &\to \Dbceent[1,-d_X]{Y},\label{eq.Rf*2}\\
\Dbceentinv[1,0]{X} &\to
\Dbceentinv[1,0]{Y},\label{eq.Rf*3}\\ 
\Dbcentinv{X} &\to \Dbcentinv[d_X]{Y}.\label{eq.Rf*4}
\end{align}
Sans hypothèse de séparation de $f$, \eqref{eq.Rf*1},
\eqref{eq.Rf*3} et \eqref{eq.Rf*4} sont encore vrais.
\end{theo}

L'hypothèse de séparation est également superflue pour
\eqref{eq.Rf*2}. On peut l'éliminer ou bien en étudiant la
$q$-divisibilité en dehors d'un sous-schéma de dimension fixée,
ou bien en utilisant une théorie de $Rf_!$ sans hypothèse de
séparation (voir \ref{coro.Rf*2}).

\begin{theo}\label{theo.Rf!h}
Soient $f: X\to Y$ un morphisme \sep de \schs \tf sur $\eta$,
$d_Y=\dim Y$, $d_r = \max_{y\in |Y|} \dim f^{-1}(y)$
. Alors $Rf^!$ induit
\begin{align}
\Dbcent{Y} &\to \Dbcent[-d_r]{X},\label{eq.Rf!1}\\
\Dbceent[1,0]{Y} &\to \Dbceent[1,-d_Y]{X},\label{eq.Rf!2}\\
\Dbceentinv[1,0]{Y} &\to \Dbceentinv[1,d_r]{X},\label{eq.Rf!3}\\
\Dbcentinv{Y} &\to \Dbcentinv[d_Y]{X}.\label{eq.Rf!4}
\end{align}
\end{theo}

Soient $X$ un \sch \tf sur $\eta$, $a_X : X\to \eta$. Rappelons que
$Ra_X^!\Qlb$ est globalement défini (pas de problème dans le cas
séparé, dans le cas général par \cite[3.2.4]{BBD}). On pose $D_X
= R\cHom(-,Ra_X^!\Qlb)$.

\begin{theo}\label{theo.D}
Soient $X$ un \sch \tf sur $\eta$, $d_X=\dim X$. Alors $D_X$ induit
\begin{align}
\Dbcentinv{X}\op &\to \Dbcent[-d_X]{X},\label{eq.D4}\\
\Dbceentinv[1,0]{X}\op &\to \Dbceent[1,0]{X},\label{eq.D3}\\
\Dbceent[1,0]{X}\op &\to \Dbceentinv[1,d_X]{X},\label{eq.D2}\\
\Dbcent{X}\op &\to \Dbcentinv{X}.\label{eq.D1}
\end{align}
De plus, pour $K\in\Modclentinv{X}$, $\fH^a(DK)$ est \ent[(a+1)],
$-d_X\le a \le -1$.
\end{theo}

\begin{theo}\label{theo.RHom}
Soient $X$ un \sch \tf sur $\eta$, $d_X=\dim X$. Alors
$R\cHom_X(-,-)$ induit
\begin{align}
\Dbcentinv{X}\op \times \Dbcent{X} &\to \Dbcent{X},\\
\Dbceentinv[1,0]{X}\op \times \Dbceent[1,0]{X} &\to \Dbceent[1,-d_X]{X},\label{eq.RHom2}\\
\Dbceent[1,0]{X}\op \times \Dbceentinv[1,0]{X} &\to \Dbceentinv[1,0]{X},\\
\Dbcentinv{X}\op \times \Dbcentinv{X} &\to
\Dbcentinv[d_X]{X}\label{eq.RHom4}.
\end{align}
\end{theo}

\section{Diviseurs \cn}\label{S3}

\begin{prop}
Soient $g: X\to Y$ un morphisme fini de \schs \tf sur $\eta$, $L\in
\Dbc(X,\Qlb)$. Alors $g_* L$ est \eent (\resp \eentinv) \ssi $L$
l'est.
\end{prop}

\begin{proof}
\Ops que $Y$ est réduit à un seul point $y$, $X$ est réduit à un
seul point $x$ et $L=\cF \in \Mod_c(X,\Qlb)$. Soient $G_y =
\Gal(\kappa(\yb)/\kappa(y))$, $G_x= \Gal(\kappa(\xb)/\kappa(x))$. Le
faisceau $\cF$ correspond à une représentation $\rho: G_x\to
\GL_{\Qlb}(\cF_\xb)$. Soient $K'$ une extension finie
quasi-galoisienne (\ie normale) de $\kappa(y)$ contenant
$\kappa(x)$, $x'=\Spec K'$. Pour $s\in G_y$, soit $\cF_s$ le
faisceau sur $x'$ correspondant à la représentation
\begin{align*}
\Gal(\kappa(\lbar{x'})/K')&\to \GL_{\Qlb}(\cF_\xb)\\
h &\mapsto \rho(s^{-1}h s).
\end{align*}
Ce faisceau ne dépend, à isomorphisme près, que de l'image de $s$
dans $G_y/G_x$. D'après la formule de Mackey (\cite[7.3]{Serre}),
on a $(g_*\cF)_{x'} \simeq \bigoplus_s \cF_s$, où $s$ parcourt un
système de représentants de $G_y/G_x$. Donc
\begin{multline*}
\text{$g_* \cF$ est $\eps(0)$-entier} \Leftrightarrow
\text{$(g_*\cF)_{x'}$ est $\eps(0)$-entier}\\
\Leftrightarrow \text{les $\cF_s$ sont $\eps(0)$-entiers}
\Leftrightarrow \text{$\cF$ est $\eps(0)$-entier.}
\end{multline*}
De même pour le cas entier inverse.
\end{proof}

Soient, en \ref{ss.purete} et \ref{lemm.deg0}, $K$ un corps
quelconque, $\eta= \Spec K$.

On va utiliser le cas spécial suivant du théorème de pureté de
Gabber \cite{Fuji}.

\begin{prop}\label{ss.purete}
Soient $n$ un entier inversible sur $\eta$, $\Lbd=\Z/n\Z$. Soit
$i:Y\to X$ une immersion fermée de \schs \regs \tf sur $\eta$
purement de codimension $c$. Alors $Ri^! \Lbd \simeq \Lbd(-c)[-2c]$.
\end{prop}

Gabber a remarqué que ce résultat découle facilement du
théorème de pureté relative \cite[XVI 3.7]{SGA4}. En effet, $i$
provient par \cb d'une immersion fermée $i_1: X_1\to Y_1$ de \schs
\tf sur $K_1$, où $K_1$ est un sous-corps de~$K$, extension \tf
d'un corps premier $K_0$. Alors $\Spec K_1$ est le point générique
d'un \sch $S_1$ intègre \tf sur $K_0$. Quitte à remplacer $S_1$
par un ouvert, \ops que $i_1$ est la fibre générique d'une
immersion fermée $i_2 : Y_2 \to X_2$ de \schs \tf sur $S_1$. Comme
$X_1$ (\resp $Y_1$) est un \sch \reg (\cite[6.5.2 (i)]{EGAIV}) et
que $X_2$ (\resp $Y_2$) est de type fini sur $K_0$, donc en
particulier, excellent, quitte à remplacer $X_2$ et $Y_2$ par des
voisinages ouverts de leurs fibres génériques, \ops que $X_2$ et
$Y_2$ sont \regs (donc lisses sur $K_0$) et $i_2$ est purement de
codimension~$c$. \Dpr le \thm de pureté relative, $Ri_2^! \Lbd
\simeq \Lbd(-c)[-2c]$. On conclut par passage à la limite.\qed

Le lemme suivant est décalqué de \cite[XXI 5.2.1]{SGA7}.

\begin{lemm}\label{lemm.deg0}
Soient $X$ un \sch \tf sur $\eta$, $a_X : X\to \eta$, $l$ un nombre
premier inversible sur $\eta$, $\cG\in \Mod_c(X,\Qlb)$.

\tmi Il existe une partie fermée $Y$ de dimension $0$ de $X$ \teq
$a_{X*} \cG \to a_{Y*}(\cG|Y)$ soit injectif, où $a_Y: Y\to \eta$.

\tmii Si $X$ est \sep de dimension~$n$, et si $U$ est un ouvert de $
X$ dont le \compl~$Z$ est de dimension $<n$, alors il existe une
partie fermée $Y$ de $U$ de dimension~$0$ et une flèche surjective
$a_{Y*} (\cG|Y)(-n) \to R^{2n}a_{X!}\cG$, où $a_Y: Y\to \eta$.
\end{lemm}

\begin{proof}
(i) est évident.

(ii) Quitte à remplacer $X$ par $X_\red$ et à rétrécir $U$, \ops
$U$ régulier purement de dimension~$n$ et $\cG|U$ lisse. Puisque
$\dim Z<n$, on a
\[0= R^{2n-1} a_{Z!} (\cG|Z) \to R^{2n}a_{U!}(\cG|U) \simto R^{2n}a_{X!}\cG
\to R^{2n}a_{Z!}(\cG|Z)=0,
\]
où $a_Z : Z\to \eta$. Donc \ops $X=U$. Appliquant (i) à
$\cGc=\cHom(\cG,\Qlb)$, on trouve une partie fermée $Y$ de $U$ de
dimension~$0$ \teq $a_{U*} \cGc \to a_{Y*} (\cGc|Y)$ soit injectif,
donc $D_\eta(a_{Y*} (\cGc|Y)) \to D_\eta(a_{U*} \cGc)$ \surj. Par le
théorème de pureté \ref{ss.purete}, on a
\[
D_\eta(a_{Y*}(\cGc|Y)) \simeq a_{Y*}(D_Y(\cGc|Y)) \simeq
a_{Y*}(\cG|Y), \] Par ailleurs, on a
\[D_\eta(a_{U*}\cGc) = \fH^0(D_\eta Ra_{U*}\cGc) \simeq
\fH^0(Ra_{U!}D_U\cGc) \simeq \fH^0(Ra_{U!} \cG(n)[2n]) =
R^{2n}a_{U!}\cG(n).\] D'où le résultat.
\end{proof}

On reprend les notations du \S~\ref{S1}.

\begin{coro}\label{coro.deg0}
Soient $X$ un \sch \tf sur $\eta$, $a_X : X\to \eta$, $\cG\in
\Mod_c(X,\Qlb)$ entier (\resp \entinv).

\tmi $a_{X*}\cG$ est \ent (\resp \entinv).

\tmii Si $X$ est \sep de dimension~$n$, alors $a_{X!}\cG$ est \ent
(\resp \entinv), $R^{2n}a_{X!}\cG$ est \ent[n] (\resp \entinv[n]).
\end{coro}

La proposition suivante est décalqué de \cite[XXI 5.3 (a)]{SGA7}.

\begin{prop}\label{prop.j*}
Soient $j: X\hra Y$ une immersion ouverte de \schs \tf sur~$\eta$ de
dimension~$1$, $\cG\in \Mod_c(X,\Qlb)$ entier. Alors $j_*\cG$ est
entier.
\end{prop}

\begin{proof}
On se ramène au cas $Y$~affine, puis $Y$~projectif.

Définissons $\fH$ par la suite exacte courte
\[0\to j_! \cG \to j_* \cG \to \fH \to 0,\]
d'où la suite exacte
\[a_{X*}\cG \to a_{Y*}\fH \to R^1a_{X!} \cG, \]
où $a_X : X\to \eta$, $a_Y: Y\to \eta$. 
D'après \ref{coro.deg0} (i), $a_{X*}\cG$ est \ent. D'après le
théorème de Deligne-Esnault \eqref{eq.Rf!b1}, $R^1 a_{X!}\cG$ est
\ent. Donc $a_{Y*}\fH$ l'est aussi. Mais $\fH$ est à supports dans
une partie fermée de dimension~$0$ de~$Y$, donc $\fH$ est \ent.
\end{proof}

\begin{prop}\label{prop.Rj*}
Soient $X$ un \sch \reg \tf sur $\eta$ de dimension~$1$, $D$ un
diviseur positif \reg. Posons $U=X-D$, $j: U\hookrightarrow X$. Soit
$\cG$ un $\Qlb$-\faisc lisse sur $U$, entier, \mr sur $X$. Alors
$Rj_* \cG$ est \eent[1,0].
\end{prop}

\begin{proof}
La question est locale sur $X$. Soit $x\in |D|$. Montrons que $(Rj_*
\cG)_x$ est \eent[1,0].

On a $\cG\simeq (\cG_\cO \ot_\cO E)\ot_E \Qlb$, où $E$ est un corps
extension finie de $\Q_l$, $\cO$ son anneau des entiers, $\cG_\cO$
un $\cO$-\faisc lisse (\constr) sur $U$. En vertu du lemme
d'Abhyankar \cite[XIII 5.2]{SGA1}, il existe, au voisinage de $x$,
un revêtement fini $g: \Xt \to X$ de la forme $\Xt=X[T]/(T^n-t)$
où $t$ est une équation locale de $x$, $n$ est un entier premier
à l'exposant caractéristique de $K$, tel que
$(g|U)^*(\cG_\cO\ot_\cO (\cO/l^2\cO))$ se prolonge en un \faisc
\loct \cst sur~$\Xt$. Comme $\cG$ est facteur direct de $(g|U)_*
(g|U)^* \cG$, on est ramené à montrer le lemme pour le faisceau
$(g|U)_* (g|U)^* \cG$. Comme $g^{-1}(D)_\red$ est un diviseur
régulier, on peut alors se ramener à montrer le lemme pour un
faisceau $\cG$ tel que $\cG_\cO \ot_\cO (\cO/l^2\cO)$ se prolonge en
un \faisc \loccst sur $X$, puis au cas $\cG_\cO \ot_\cO
(\cO/l^2\cO)$ constant par la formule de projection.

Soient $X_\xp$ le hensélisé de $X$ en $x$, $U_\xp = X_\xp \times_X
U$, $j_\xp: U_\xp \hookrightarrow X_\xp$, $\fH = \cG|U_\xp$. Alors
$\fH \simeq (\fH_\cO\ot_\cO E)\ot_E \Qlb$, avec $\fH_\cO\ot_\cO
(\cO/l^2\cO)$ \cst. On a $(Rj_{\xp*} \fH)_x = (Rj_* \cG)_x \in
\Dbc(x, \Qlb)$. D'après \ref{prop.j*}, $(j_{\xp*}\fH)_x =
(j_*\cG)_x$ est \ent. Il reste à montrer que $(R^1 j_{\xp*}\fH)_x$
est \ent[1].

On a une suite exacte de groupes
\[1\to \Zch(1) \to G \to \Gal(\kappa(\xb)/\kappa(x)) \to 1,\]
où $\Zch(1) = \prod_{p'\neq \car(K)} \Z_{p'}(1)$, $G=
\pimod(U_\xp)$. Le \faisc $\fH$ correspond à une \repnl de $G$.
D'après le \thm de monodromie locale, la restriction de cette \repn
à $\Zch(1)$ est quasi-\unipe, donc unipotente. On obtient une
filtration $M$ finie, croissante de $\fH$, \teq chaque $\gr^M_a\fH$
se prolonge en un $\Qlb$-\faisc lisse $\cG_a$ sur $X_\xp$.

Montrons que $(R^1j_{\xp*} M_a )_x$ est \ent[1] par \rec sur $a$, ce
qui achèvera la \dem de la proposition. L'assertion est claire pour
$a\ll 0$. Supposons l'assertion établie pour $a-1$. La suite exacte
courte
\[0 \to M_{a-1} \to M_a \to \cG_a|U_\xp \to 0\]
donne le triangle \dist
\[Rj_{\xp*} M_{a-1} \to Rj_{\xp*} M_a  \to Rj_{\xp*}(\cG_a|U_\xp)
\to.\]
D'après une formule de projection, $Rj_{\xp*}(\cG_a|U_\xp) \simeq
Rj_{\xp*}\Qlb \otimes \cG_a$. On a
\[(R^qj_{\xp*}\Qlb)_x =
\begin{cases}
\Qlb(-q) & \text{si $q=0,1$,}\\
0 & \text{sinon.}
\end{cases}
\] Donc on a la suite exacte
\[(j_{\xp*}M_a )_x \to (\cG_a)_x \to (R^1 j_{\xp*}M_{a-1})_x
\to (R^1 j_{\xp*}M_a )_x \to (\cG_a)_x (-1).\] Ici $(j_{\xp*}M_a
)_x$ est un sous-\faisc de $(j_{\xp*}\fH )_x$, donc \ent. Par \hyp
de \rec, $(R^1 j_{\xp *} M_{a-1})_x$ est \ent[1]. Donc $(\cG_a)_x$
est entier, $(\cG_a)_x(-1)$ est $1$-entier. On en déduit que $(R^1
j_{\xp *} M_{a})_x$ est \ent[1].
\end{proof}

Le lemme suivant est une variante de \cite[XXI 5.6.2]{SGA7}.

\begin{lemm}\label{lemm.cb}
Soient $X$ un \sch \noeth \reg, $D=\sum_\idI D_i$ un \dscn de $X$
avec $(D_i)_\idI$ une famille finie de diviseurs \regs,  $U=X-D$,
$n$ un entier inversible sur $X$, $\Lbd=\Z/n\Z$, $\cG\in
\Mod_c(U,\Lbd)$ \loct\cst, \mr sur $X$.

\tmi Soient $i\in I$, $U_\ip =X-\bigcup_{h\in I-\{i\}}D_h$,
$D_{i,U_\ip} = D_i \times_X U_\ip$, d'où un diagramme à carré
cartésien
\[\xymatrix{&D_{i,U_\ip}\ar@{^{(}->}[r]^{j'_{\ip}}\ar[d]^{\iota'_{i}} &
D_i\ar[d]^{\ii}\\
U\ar@{^{(}->}[r]^{j^\ip} & U_\ip\ar@{^{(}->}[r]^{j_\ip}& X}\]
Alors le morphisme de \cb
\beq\label{eq.lemm.cb}
\ii^* R j_{\ip*}(R j^\ip_{*}\cG) \to Rj'_{\ip*}
\iota'^*_{i}(Rj^\ip_{*} \cG)
\eeq
est un \isom et les \faisx
\[\iota'^*_{i} R^qj^\ip_{*} \cG, \quad q\in \Z,\]
sont \loct\csts, \mrs sur $D_i$.

\tmii Soit $f:Y\to X$ un \morph de \schs \regs \noeths. Supposons
que $f^{-1}(D)$ soit un \dcn et que $(f^{-1}(D_i))_{i\in I}$ soit
une famille de diviseurs réguliers. Considérons le \cacart :
\[\xymatrix{Y_U \ar@{^{(}->}[r]^{j_Y} \ar[d]^{f_U} &
Y\ar[d]^{f}\\
U\ar@{^{(}->}[r]^{j} & X}\] Alors le \morph de \cb $f^* Rj_* \cG
\to Rj_{Y*} f^*_U \cG$ est un \isom.
\end{lemm}


\begin{proof}
La question est locale sur $X$. Soit $x$ un point de $X$. En vertu
du lemme d'Abhyankar, il existe, au voisinage de $x$, un revêtement
fini $g: \Xt \to X$ de la forme
\[\Xt=X[T_1,
\dots,T_r]/(T_1^{n_1}-t_1, \dots, T_r^{n_r}-t_r)\]
où les $t_i$
sont des équations locales des composantes de $D$ passant par $x$,
et $n_i$ des entiers premiers à l'exposant caractéristique de
$\kappa(x)$, \tq $(g|U)^*\cG$ se prolonge en un \faisc \loct \cst
sur $\Xt$. Comme $\cG$ s'injecte dans $(g|U)_* (g|U)^* \cG$ et le
quotient $\cG_1$ est \mr sur $X$, on peut itérer cette
construction. Pour tout $N\ge 1$, on obtient, quitte à rétrécir
$X$, une résolution
\[\cG \to  (g|U)_* (g|U)^*
\cG \to (g_1|U)_* (g_1|U)^* \cG_1 \to \dots \to (g_N|U)_* (g_N|U)^*
\cG_N.
\]
Donc on est ramené à montrer le lemme pour le faisceau $(g|U)_*
(g|U)^* \cG$. Comme $g^{-1}(D)_\red = \sum_{i\in I}
g^{-1}(D_i)_\red$ est un \dcn avec $(g^{-1}(D_i)_\red)_{i\in I}$ une
famille de diviseurs \regs, on peut alors se ramener à montrer le
lemme pour un faisceau $\cG$ qui se prolonge en un \faisc \loccst
sur $X$, puis au cas $\cG=\Lbd_U$ par la formule de projection.

Le point (ii) résulte alors de \cite[\S~8]{Fuji} et de la
fonctorialité des classes des diviseurs \cite[Th.\ finitude,
2.1.1]{SGA4d}.

Pour (i), notons que $D_{i,U_\ip}$ est un diviseur \reg de $U_\ip$,
de complémentaire~$U$. Pour tout~$q$,
\[\iota'^*_i R^q j^\ip_* \Lbd_U \simeq
\begin{cases}
\Lbd_{D_{i,U_\ip}} &\text{si $q=0$,}\\
\Lbd_{D_{i,U_\ip}}(-1) &\text{si $q=1$,}\\
0&\text{sinon.}
\end{cases}\] est \loct\cst, \mr sur $D_i$. On a donc un triangle \dist
\[\iota'_{i*} \Lbd_{D_i,{U_\ip}}(-1)[-2] \to \Lbd_{U_\ip} \to Rj^\ip_* \Lbd_U \to. \]
Le morphisme \eqref{eq.lemm.cb} est un \isom car le morphisme de
changement de base $\ii^* Rj_{\ip *} \to Rj'_{\ip *} \iota'^*_{i}$
induit des \isoms sur $\Lbd_{U_\ip}$ en vertu de (ii) et
trivialement sur $\iota'_{i*} \Lbd_{D_i,{U_\ip}}$.
\end{proof}

La proposition suivante est décalquée de \cite[XXI 5.6.1]{SGA7}.

\begin{prop}\label{prop.Rj*dcn}
Soient $X$ un \sch \tf sur $\eta$, $D$ un \dcn. Posons $U=X-D$,
$j: U\hookrightarrow X$. Soit $\cG$ un $\Qlb$-\faisc lisse sur
$U$, \ent, \mr sur $X$. Alors $Rj_*\cG$ est \eent[1,0].
\end{prop}

\begin{proof}
Le problème étant local pour la \topo étale au voisinage d'un
point fermé de $D$, \ops $D$ \scn. Comme $\Reg(X)$ est un ouvert de
$X$ 
contenant $D$, \ops $X$ \reg.

On pose $D=\sum_\idI D_i$ avec $(D_i)_\idI$ une famille finie de
diviseurs \regs. On fait une \rec sur $n=\sharp I$. Le cas $n=0$ est
trivial. Pour $n>0$, on choisit $i\in I$ et applique \ref{lemm.cb}
(i), dont on conserve les notations. Pour tout $x\in
\left|\DiUi\right|$, il existe un sous-schéma \reg de $U_\ip$ de
dimension~$1$ \tq son intersection avec $\DiUi$ soit le schéma $x$.
Donc $Rj^\ip_* \cG$ est \eent[1,0], d'après \ref{lemm.cb} (ii) et
\ref{prop.Rj*}. Notons que $\sum_{h\in I-\{i\}} D_h\cap D_i$ est un
diviseur à croisements normaux de $D_i$ de complémentaire
$D_{i,U_\ip}$, et pour tout~$q$, $\iota'^*_i R^q j^\ip_* \cG$ est un
$\Qlb$-\faisx lisse sur $D_{i,U_\ip}$, \mr sur $D_i$ en vertu de
\ref{lemm.cb} (i). Donc $\iota_i^* Rj_* \cG \simeq Rj'_{\ip*}
\iota'^*_i Rj^\ip_* \cG$ est \eent[1,0], d'après l'\hyp de \rec.
Comme $i$~est arbitraire, on en conclut que $Rj_* \cG$ est
\eent[1,0].
\end{proof}

\section{Démonstration de \ref{theo.Rf!b} à \ref{theo.RHom}}\label{S4}

La proposition suivante est une variante de \cite[2.6]{Orgogozo}.

\begin{prop}\label{lemm.dJ}
Soient $F$ un corps, $X$ un \sch \sep \tf sur $\Spec F$, $U$ une
partie ouverte 
de $X$.

\tmi Il existe un morphisme $r_0 : X'_0 \to X$ propre surjectif
avec $X'_0$ \reg 
et un sous-\sch ouvert fermé $W_0$ de
$X'_0$ contenant $r_0^{-1}(U)$ 
\tsq $r_0^{-1}(U)$ soit le \compl
d'un \dscn 
dans $W_0$
.

\tmii Pour tout $n\ge 0$, il existe une \ext finie radicielle $F'$
de $F$ et un hyperrecouvrement propre \ntron \sscin $r\sm :X'\sm \to
X_{F'}$ \tsq $X'_m$ soit lisse sur $\Spec F'$ et que
$r_m^{-1}(U_{F'})$ soit le
\compl d'un \dscn \relta $\Spec F'$ 
dans un
\ssch ouvert fermé de $X'_m$
, $0\le m \le n$.
\end{prop}

\begin{proof}
(i) Au cas où $X$ est intègre et $U\neq\emptyset$, il existe un
morphisme $r_0 : X'_0 \to X$ propre surjectif avec $X'_0$ intègre et
\reg \tq $r_0^{-1}(U)$ soit le \compl d'un \dscn, en vertu de
\cite[4.1]{deJong}. On prend $W_0 = X'_0$.

Le cas où $X$ est intègre et $U=\emptyset$ en résulte : appliquer
le cas précédent à $X$ et la partie ouverte~$X$ pour obtenir
$r_0$, et puis prendre $W_0 = \emptyset$.

Dans le cas \gel, soient $X_\al$ les schémas réduits associés aux
composantes \irrs de $X$, $a: \coprod X_\al \to X$ le morphisme
\can. Alors $a$ est fini et \surj. Pour chaque $\al$, appliquons (i)
à $X_\al$ et $U\times_X X_\al$, on obtient $\phi_\al : (X_\al)'_0
\to X_\al$ propre \surj avec $(X_\al)'_0$ \reg et un sous-\sch
ouvert fermé $W_\al$ de $(X_\al)'_0$ contenant l'image inverse
$U_\al$ de $U$ \tsq $U_\al$ soit le complémentaire dans $W_\al$
d'un \dscn. Posons $X'_0 = \coprod(X_\al)'_0$, $r_0 = a\circ
\coprod\phi_\al$, $W_0 = \coprod W_\al$. Alors $r_0$ et $W_0$
satisfont aux conditions de (i).

(ii) Cas $F$ parfait. On fait une \rec sur $n$. Lorsque $n=0$, (ii)
dégénère en (i). Supposons donné un \hyperr \ntron $r\sm : X'\sm
\to X$ \ver ant les conditions de (ii). On applique (i) au $X$-\sch
$(\cosq_n X'\sm)_{n+1}$ et l'image inverse de $U$
. On obtient $\beta : N\to (\cosq_n X'\sm)_{n+1}$ propre \surj
avec $N$ lisse sur $\Spec
F$ 
et un \ssch ouvert fermé $W$ de $N$ 
\tsq l'image inverse de $U$ dans $N$ soit le \compl d'un \dscn
\relta $\Spec F$ dans 
$W$
. L'\hyperr propre $(n+1)$-tronqué \sscin associé au triplet
$(X'\sm,N,\beta)$ \cite[V\up{bis} 5.1.3]{SGA4} \ver e les
conditions de (ii) pour $n+1$.

Cas \gel. On prend une \clot parfaite $\Fb$ de $F$ et applique
(ii) à $\Fb$, $X_{\Fb}$ et $U_{\Fb}$. L'\hyperr tronqué et les
diviseurs \scn obtenus se descendent à une sous-\ext finie $F'$ de
$F$.
\end{proof}

\begin{ssect}\label{ss.demoRf*}
\begin{proof}[Démonstration de \eqref{eq.Rf*1}]
Il faut montrer que pour un $\Qlb$-\faisc \constr $\cG$ sur $X$,
\ent, $Rf_* \cG$ est \ent.

On fait une \rec sur $d=\dim X$. Le cas $d\le 0$ est trivial.

Soit $d\ge 1$. 
Choisissons un ouvert
affine $U\stackrel{j}{\hookrightarrow} X$ 
\tq $\cG|U$ soit lisse et que son \compl
$Z\stackrel{i}{\hookrightarrow} X$ soit de dimension $<d$. Le
triangle distingué
\[
i_* Ri^! \cG \to \cG \to Rj_* j^* \cG \to
\]
induit le triangle distingué
\[R(fi)_* Ri^! \cG \to Rf_* \cG \to
R(fj)_* j^* \cG \to.
\]
Compte tenu de l'\hyp de \rec, il suffit de voir que $Rj_* j^* \cG$
et $R(fj)_* j^* \cG$ sont \ents. Il suffit donc de vérifier le \thm
sous l'\hyp
additionnelle que $X$ est \sep 
et $\cG$ lisse.

On a $\cG\simeq (\cG_\cO \ot_\cO E)\ot_E \Qlb$ avec $\cG_\cO$ lisse.
Soit $p: X'\to X$ un \rev \etale surjectif qui trivialise
$\cG_\cO\ot_\cO (\cO/\fm)$, où $\fm$ est l'idéal maximal de $\cO$.
Le faisceau $\cG$ est facteur direct de $p_* p^* \cG$, de sorte
qu'il suffit de voir l'\inte de $R(fp)_* p^* \cG$. Donc il suffit de
\ver er le \thm sous l'\hyp additionnelle que $X$ est \sep et
\begin{numcond}\label{cond.theo}
$\cG\simeq (\cG_\cO \ot_\cO E)\ot_E \Qlb$ avec $\cG_\cO\ot_\cO
(\cO/\fm)$ constant.
\end{numcond}

On factorise $f$ en $X\hookto{j} Z\sto{g} Y$, où $j$ est une
immersion ouverte
, $g$ est un morphisme propre. Comme $Rg_*$ préserve
l'intégralité en vertu du théorème de Deligne-Esnault
\eqref{eq.Rf!b1}, il suffit de prouver l'intégralité de $Rj_*\cG$.
On est donc ramené à démontrer \eqref{eq.Rf*1} pour $j$ et $\cG$.
Pour cela, \ops $Z$ affine, donc \sep.

Soit $i\ge 0$. On applique \ref{lemm.dJ} (ii) à $j$ et $n=i+1$.
Quitte à changer les notations, on peut supposer que l'extension
radicielle de \emph{loc.\ cit.}\ est triviale. On obtient un carré
cartésien (de \schs simpliciaux $(i+1)$-tronqués)
\[\xymatrix{X'\sm \ar@{^{(}->}[r]^{j'\sm}\ar[d]^{s\sm} &Z'\sm\ar[d]^{r\sm}\\
X\ar@{^{(}->}[r]^{j} &Z}\] où $r\sm$ est un \hyperr propre
$(i+1)$-tronqué \sscin, $Z'_m$ lisse sur~$\eta$, $j'_m$ une
immersion ouverte faisant de $X'_m$ le \compl d'un \dcn \relta
$\eta$ dans une partie ouverte fermée de $Z'_m$, $0 \le m \le i+1$.
D'après la descente cohomologique,
\[\tau_{\le i} Rj_* \cG \simeq  \tau_{\le i} R j_* R{s\sm}_* s\sm^*
\cG = \tau_{\le i} R {r\sm}_* R{j'\sm}_* s\sm^* \cG.
\]
Comme $R{r\sm}_*$ préserve l'intégralité \eqref{eq.Rf!b1}, il
suffit de voir l'\inte des $Rj'_{m*}s_m^* \cG$, $0\le m \le i$. Or
$s^*_m \cG$ satisfait encore à \ref{cond.theo}, donc est \mr sur
$Z'_m$. Il suffit alors d'appliquer \ref{prop.Rj*dcn}.
\end{proof}
\end{ssect}

\begin{prop}\label{prop.DT!}
Soient $X$ un \sch \reg \sep \tf sur $\eta$, purement de
dimension~$1$, $a_X: X\to \eta$, $\cG$ un $\Qlb$-\faisc lisse sur
$X$, \entinv. Alors $R^1 a_{X!}\cG$ est $1$-entier inverse.
\end{prop}

\begin{proof}
\Ops que $\cG \simeq (\cG_\cO\otimes_\cO E)\otimes_E \Qlb$, avec
$\cG_\cO$ lisse. Pour $p : X'\to X$ un revêtement étale \surj qui
trivialise $\cG_\cO \ot_\cO (\cO/\fm)$, $\cG$ est un facteur direct
de $p_*p^*\cG$. Ceci permet de supposer $\cG_\cO \otimes_\cO
\cO/\fm$ constant. On a
\[D_\eta(R^1a_{X!}\cG) \simeq \fH^{-1}(D_\eta Ra_{X!}\cG) \simeq
\fH^{-1}(Ra_{X*}\cGc(1)[2]) = R^1a_{X*}\cGc(1).\] Soit $j: X\hra P$
une compactification régulière de $X$. De $a_X = a_P \circ j$ on
déduit une suite spectrale
\[E_2^{pq} = R^pa_{P*}R^qj_*\cGc \Rightarrow R^{p+q} a_{X*} \cGc.\]
D'après \ref{prop.Rj*dcn} et \ref{coro.deg0}, $R^0
a_{P*}R^1j_*\cGc$ est \ent. D'autre part, $R^1 a_{P*}R^0j_*\cGc$ est
un quotient de $R^1a_{P*} j_! \cGc$, donc entier en vertu de
\eqref{eq.Rf!b1}. Donc $R^1a_{X*}\cGc$ est entier, $R^1a_{X!}\cG$
est $1$-entier inverse.
\end{proof}

\begin{proof}[Démonstration de \ref{theo.Rf!b}]
Comme on a déjà remarqué, le cas \og entier \fg (\eqref{eq.Rf!b1}
et \eqref{eq.Rf!b2}) de \ref{theo.Rf!b} est démontré dans
\cite[0.2]{DE}. Supposons maintenant $\cF$ entier inverse. \Ops
$Y=\eta$, $f=a_X$.

Traitons d'abord le cas $X=\bA^1_{\eta}$. Soient $j:
U\hookrightarrow X$ un ouvert dense \tq $\cF|U$ soit lisse, $i:
Z\hra X$ le fermé complémentaire. La suite exacte
\[0\to j_! j^* \cF \to \cF \to i_* i^* \cF \to 0\]
donne le triangle \dist
\[Ra_{U!} (\cF|U) \to Ra_{X!}\cF \to Ra_{Z!}(\cF|Z) \to .\]
D'après \ref{coro.deg0}, $Ra_{Z!} (\cF |Z)$ est \entinv, $R^0 a_{U!}
(\cF|U)$ est \entinv, $R^2 a_{U!} (\cF|U)$ est \entinv[1]. D'après
\ref{prop.DT!}, $R^1 a_{U!} (\cF|U)$ est \entinv[1]. Donc $Ra_{X!}
\cF$ est \entinv[I] et \entinv[1].

Pour le cas général, procédons par \rec sur $n$. Le cas $n\le 0$
est trivial. Soit $n\ge 1$. D'après le lemme de normalisation, il
existe $j: U\hookrightarrow X$ un ouvert dense de $X$ et un
morphisme $f: U\to Y=\bA^1_\eta$ à fibres de dimension $\le n-1$.
Soit
 $i: Z\hra X$ le \compl de $U$. Alors $Z$ est de dimension $\le n-1$. La suite exacte
\[0\to j_! j^* \cF \to \cF \to i_* i^* \cF \to 0\]
donne le triangle \dist
\[Ra_{U!} (\cF|U) \to Ra_{X!}\cF \to Ra_{Z!}(\cF|Z) \to .\]
Par l'\hyp de \rec, $Ra_{Z!}(\cF|Z)$ est \eentinv[1,0] et
\entinv[(n-1)]. Il suffit donc de \ver er la \propo pour $\cF|U$.
Ceci résulte de la suite spectrale
\[E_2^{pq}= R^p
a_{Y!} R^q f_! (\cF|U) \Rightarrow R^{p+q}a_{U_!} (\cF|U),\] de
l'\hyp de \rec appliquée aux fibres de $f$ et du cas d'une droite
affine déjà traité.
\end{proof}

La proposition suivante est un analogue de \ref{prop.Rj*}.

\begin{prop}\label{prop.Rj*inv}
Soient $X$ un \sch \reg \tf sur $\eta$ de dimension~$1$, $D$ un
diviseur positif \reg. Posons $U=X-D$, $j: U\hookrightarrow X$. Soit
$\cG$ un $\Qlb$-\faisc lisse sur $U$, entier inverse, \mr sur $X$.
Alors $Rj_* \cG$ est \eentinv[1,0].
\end{prop}

\begin{proof}
Comme $\cGc$ est lisse, \mr sur $X$, \ent, $Rj_* \cGc$ est
\eent[1,0] en vertu de \ref{prop.Rj*dcn}. Soit $i : D\hookrightarrow
X$. D'après la dualité locale en dimension~$1$, $i^* R^1 j_* \cG
\simeq (i^* j_* \cGc)\spcheck (-1)$ est $1$-entier inverse, $i^*
j_*\cG \simeq (i^* R^1 j_* \cGc)\spcheck (-1) $ est \entinv.
\end{proof}

La proposition suivante est un analogue de \ref{prop.Rj*dcn}.

\begin{prop}\label{prop.Rj*dcninv}
Soient $X$ un \sch \tf sur $\eta$, $D$ un \dcn. Posons $U=X-D$,
$j: U\hookrightarrow X$. Soit $\cG$ un $\Qlb$-\faisc lisse sur
$U$, \entinv, \mr sur $X$. Alors $Rj_*\cG$ est \eentinv[1,0].
\end{prop}

On déduit \ref{prop.Rj*dcninv} de \ref{prop.Rj*inv}  de la \me
manière qu'on a déduit \ref{prop.Rj*dcn} de \ref{prop.Rj*}.

On déduit \eqref{eq.Rf*3} de \ref{prop.Rj*dcninv} de la \me manière
qu'on a déduit \eqref{eq.Rf*1} de \ref{prop.Rj*dcn} en
\ref{ss.demoRf*}.

\begin{remq}\label{remq.Ri!}
L'assertion \eqref{eq.Rf!1} (\resp \eqref{eq.Rf!3}) pour $f$ une
immersion fermée découle de ce qui précède. En effet, soient $j: Y-X
\hra Y$ l'ouvert complémentaire, $K\in\Dbc(Y,\Qlb)$ \ent (\resp
\entinv[I]). On a le \trdist
\[Rf^! K \to f^*K \to f^* Rj_* j^* K \to.\]
En appliquant \eqref{eq.Rf*1} (\resp \eqref{eq.Rf*3}) à $j$, on
obtient que $f^* Rj_* j^* K$ est \ent (\resp \entinv[I]), donc $f^*
Rj_* j^* K [-1]$ l'est aussi (\resp \entinv[(I-1)]). Or $f^*K$ est
\ent (\resp \entinv[I]). On en conclut que $Rf^! K$ l'est aussi.
Cette démonstration donne un peu plus dans le cas entier inverse :
pour $K\in \Modclentinv{Y}$, $R^a f^! K$ est \entinv[(a-1)], $a\ge
1$.
\end{remq}

\begin{ssect}\label{ss.demo.D}
\begin{proof}[Démonstration de \ref{theo.D}]
\Ops $X$ réduit. On fait une \rec  sur $d_X$. Le cas $d_X \le 0$
est trivial. Pour $d_X \ge 1$, il suffit de montrer que pour
$K\in\Mod_c(X,\Qlb)$ entier (\resp \entinv), $DK$ est \entinv et
\eentinv[1,d_X] (\resp \eent[1,0] et \ent[-d_X] et $\fH^a(DK)$ est
\ent[(a+1)], $-d_X\le a \le -1$). Prenons un ouvert $j: U\hra X$
\reg purement de dimension $d_X$ \tq le \compl $i: V\hra X$ soit de
dimension $<d_X$ et que $K|U$ soit lisse. Alors $D(j^*K)\simeq (j^*
K)\spcheck(d_X)[2d_X]$.

On a le triangle \dist
\beq\label{eq.demo.D}
i_* D(i^*K) \to DK \to Rj_*D(j^*K) \to .
\eeq
\Dpr \eqref{eq.Rf*3}
(\resp \eqref{eq.Rf*1}), $Rj_*D(j^*K)$ est \entinv[(I+d_X)] (\resp
\ent[-d_X]). \Dpr l'\hyp de \rec, $i_* D(i^*K)$ l'est aussi. Donc
$DK$ est \eentinv[1,d_X] (\resp \ent[-d_X]).

On a le triangle \dist
\[j_! D(j^*K) \to DK \to i_*D(Ri^!K) \to. \]
Le terme $j_! D(j^*K)\in D^{[-2d_X,-2d_X]}$ est \entinv (\resp
\ent[I]). D'après \ref{remq.Ri!}, $Ri^!K$ est \ent (\resp $R^0 i^!
K$ est \entinv et $R^a i^! K$ est \entinv[(a-1)], $a\ge 1$). \Dpr
l'\hyp de \rec, $D(Ri^!K)$ est \entinv (\resp $\fH^0(D(Ri^!K))$ est
\ent et $\fH^a(D(Ri^!K))$ est \ent[(a+1)], $a\le -1$, en vertu de la
suite spectrale \[E_2^{pq} = \fH^p(D(R^{-q}i^!K)) \Rightarrow
\fH^{p+q}(DRi^!K)).\] Donc $DK$ est \entinv (\resp \eent[1,0] et
$\fH^a(DK)$ est \ent[(a+1)], $-d_X\le a \le -1$).
\end{proof}
\end{ssect}

\begin{proof}[Démonstration de \ref{theo.Rf*b}]
Les assertions \eqref{eq.Rf*1}
et \eqref{eq.Rf*3} sont déjà démontrées en \ref{ss.demoRf*} et \ref{prop.Rj*dcninv}.
Si $f$ est \sep, alors \eqref{eq.Rf*2} et \eqref{eq.Rf*4} découlent
de \ref{theo.Rf!b} et \ref{theo.D} : $Rf_* \simeq D_Y Rf_! D_X$
induit
\begin{gather*}
\Dbceent[1,0]{X} \sto{D_X} \Dbceentinv[1,d_X]{X}\op \sto{Rf_!}
\Dbceentinv[1,d_X]{Y}\op \sto{D_Y}
\Dbceent[1,-d_X]{Y},\\
\Dbcentinv{X} \sto{D_X} \Dbcent[-d_X]{X}\op \sto{Rf_!}
\Dbcent[-d_X]{Y}\op \sto{D_Y} \Dbcentinv[d_X]{Y}.
\end{gather*}
Il reste à montrer \eqref{eq.Rf*4} sans supposer $f$ \sep. Pour
cela, prenons un recouvrement ouvert affine $W\to X$. Soit $g\sm:
\cosq_0(W/X)\sm \to X$. Alors $Rf_* \simeq Rf_* R{g\sm}_* g\sm^*$,
et il suffit d'appliquer le résultat du cas \sep.
\end{proof}

\begin{remqe}
Les assertions \eqref{eq.Rf!2} et \eqref{eq.Rf!4}
découlent 
de \ref{theo.D} : $Rf^! \simeq D_X f^* D_Y$ induit
\begin{gather*}
\Dbceent[1,0]{Y} \sto{D_Y} \Dbceentinv[1,d_Y]{Y}\op \sto{f^*}
\Dbceentinv[1,d_Y]{X}\op \sto{D_X}
\Dbceent[1,-d_Y]{X},\\
\Dbcentinv{Y} \sto{D_Y} \Dbcent[-d_Y]{Y}\op \sto{f^*}
\Dbcent[-d_Y]{X}\op \sto{D_X} \Dbcentinv[d_Y]{X}.
\end{gather*}
\end{remqe}

Dans le cas $f$ quasi-fini, on peut améliorer \eqref{eq.Rf*2}
comme suit.

\begin{prop}\label{prop.Rf*qf}
Soit $f: X\to Y$ un morphisme quasi-fini séparé de \schs \tf sur
$\eta$ avec $d_X=\dim X\ge 1$. Alors $Rf_*$ envoie
$\Dbcleent[1,0]{X}$ dans $\Dbcleent[1,1-d_X]{Y}$.
\end{prop}

\begin{proof}
Par le théorème principal de Zariski, \ops que $f$ est une
immersion ouverte dominante. Soit $K\in \Mod_c(X,\Qlb)$ \ent. On a
le \trdist
\[i_* Ri^! f_! K \to f_! K \to Rf_*K \to,\]
où $i: Y-X\hra Y$ est le fermé complémentaire. Il suffit donc
d'appliquer \eqref{eq.Rf!2}.
\end{proof}

\vcb{v1m}{La proposition suivante généralise \cite[0.4]{DE}.}

\begin{prop}\label{prop.i^!}
Soient $i: X\hra Y$ une immersion de \schs \tf sur~$\eta$ avec $Y$
\reg, $d_Y = \dim Y $, $d_c=\codim(X,Y)$, $\cG\in \Mod_c(Y,\Qlb)$
entier (\resp \entinv) lisse. Soit $\eps: \Z \to \Q$ une fonction
\ver ant
\[\eps(a)=
\begin{cases}
d_c & \text{si $2d_c \le a<d_c+d_Y$,}\\
a+1-d_Y & \text{si $d_c+d_Y \le a < 2d_Y$,}\\
d_Y & \text{si $a=2d_Y$.}\\
\end{cases}\]
Alors $Ri^!\cG$ est \ent[\eps]
(\resp \entinv[(I-d_c)]).
\end{prop}

\begin{proof} 
On va montrer que pour toute partie fermée $W$ de $X$, $Ri_W^!
\cG$ est \ent[\eps] (\resp \entinv[(I-d_c)]), où $i_W : W\hra Y$.
Ici on a muni $W$ de la structure de \sch réduit induite. On fait
une \rec noethérienne. Le cas $W=\emptyset$ est trivial. Pour
$W\neq \emptyset$, on prend un ouvert \reg \irr $j: U\hra W$. Soit
 $i_Z : Z\hra W$ son \compl. On a le triangle
\dist
\[i_{Z*}Ri_Z^! Ri_W^! \cG \to Ri_W^!\cG \to Rj_* j^* Ri_W^!\cG \to.\]
D'après \ref{ss.purete}, $j^* Ri_W^!\cG = R(i_Wj)^!\cG
\simeq \cG(-d)[-2d]$, où $d=\codim(U,Y)\ge d_c$.  
Si $d_U = \dim U=0$, alors $Rj_* j^* Ri_W^!\cG$ est \eent car
$\eps(2d)\le d$ ; si $d_U \ge 1$, alors \dpr \eqref{eq.Rf*1} et
\ref{prop.Rf*qf}, $Rj_* j^* Ri_W^!\cG$ est \ent[d] et
\ent[(I+1-d-d_U)], donc \eent, compte tenu du fait que $d+d_U \le
d_Y$. (\Resp \dpr \eqref{eq.Rf*3}, $Rj_* j^* Ri_W^!\cG$ est
\entinv[(I-d)], donc \entinv[(I-d_c)].) \Dpr l'\hyp de \rec, $Ri_Z^!
Ri_W^! \cG = R(i_W i_Z)^! \cG$ est \ent[\eps] (\resp
\entinv[(I-d_c)]). Donc $Ri_W^! \cG$ l'est aussi.
\end{proof}

\begin{ssect}\label{ss.demo.Rf!h}
\begin{proof}[Démonstration de \ref{theo.Rf!h}]
Les assertions \eqref{eq.Rf!2} et \eqref{eq.Rf!4} sont traitées
plus haut. Pour le reste, on fait une \rec sur $d_Y$. Le cas $d_Y
<0$ est clair. Pour $d_Y \ge 0$, \ops $Y$ réduit. Soit $\cG\in
\Mod_c(Y,\Qlb)$ \ent (\resp \entinv). Il existe un ouvert \reg~$U$
de~$Y$ de \compl~$W$ de dimension $\le d_Y -1$ \tq $\cG|U$ soit
lisse. On considère le diagramme à carrés cartésiens
\[\xymatrix{X_U\ar[d]^{f_U} \ar@{^{(}->}[r]^{j'} & X\ar[d]^f & X_W \ar[d]^{f_W}
\ar@{_{(}->}[l]_{i'} \\
U \ar@{^{(}->}[r]^{j} & Y & W\ar@{_{(}->}[l]_{i}}\] On a le
triangle \dist
\beq\label{eq.demo.Rf!}
i'_* Rf_W^! Ri^! \cG \to Rf^!\cG \to Rj'_* Rf_U^! j^* \cG \to .
\eeq
D'après \ref{remq.Ri!}, $Ri^!\cG$ est \ent (\resp \eentinv[1,0]),
donc $Rf_W^! Ri^!\cG$ est \ent[-d_r] (\resp \eentinv[1,d_r]) en
vertu de l'\hyp de \rec. Il reste à considérer $Rf_U^!(\cG|U)$.

Il suffit de montrer que pour $Y$ \reg et $\cG\in \Mod_c(Y,\Qlb)$
\ent (\resp \entinv) lisse, $Rf^!\cG$ est \ent[-d_r] (\resp
\eentinv[1,d_r]). Le problème étant local sur $Y$, \ops $Y$ \irr et
que $f$ se factorise en $X\stackrel{i_X}{\hra} \bA_Y^n \sto{p} Y$,
où $i_X$ est une immersion fermée. $\bA_Y^n$ est \irr
\cite[4.5.8]{EGAIV}, donc \bieqdim \cite[5.2.1]{EGAIV}, donc
$\codim(X,\bA_Y^n) = n+d_Y -d_X \ge n-d_r$, car $d_X \le d_Y + d_r$.
Il suffit donc d'appliquer \ref{prop.i^!} à $i_X$ et $Rp^! \cG[-2n]
= p^*\cG(n)$.
\end{proof}
\end{ssect}

\begin{ssect}\label{ss.demo.RHom}
\begin{proof}[Démonstration de \ref{theo.RHom}]
Les assertions \eqref{eq.RHom2} et \eqref{eq.RHom4} découlent de
\ref{theo.D} :
\[R\cHom_X(-,-) \simeq D_X (-\otimes D_X -)\]
induit
\begin{gather*}
\begin{split}
\Dbceentinv[1,0]{X}\op \times \Dbceent[1,0]{X} \sto{(\id,D_X)} &
\Dbceentinv[1,0]{X}\op \times
\Dbceentinv[1,d_X]{X}\op\\
&\sto{\otimes} \Dbceentinv[1,d_X]{X}\op \sto{D_X}
\Dbceent[1,-d_X]{X},
\end{split}\\
\begin{split}
\Dbcent{X}\op \times \Dbcentinv{X} \sto{(\id,D_X)} & \Dbcent{X}\op
\times
\Dbcent[-d_X]{X}\op\\
&\sto{\otimes} \Dbcent[-d_X]{X}\op \sto{D_X} \Dbcentinv[d_X]{X}.
\end{split}
\end{gather*}
Pour le reste, il suffit de montrer que pour $K,L\in
\Mod_c(X,\Qlb)$, $K$ \entinv, $L$ \ent (\resp $K$ \ent, $L$
\entinv), $R\cHom(K,L)$ est \ent (\resp \eentinv[1,0]). Par
dévissage de $K$, on se ramène à supposer $K$ de la forme $j_!
\cG$, où $j: Y\hra X$ une immersion, $\cG\in \Mod_c(Y,\Qlb)$
\entinv (\resp \ent) lisse. Alors
\[R\cHom_X(j_!\cG,L) \simeq Rj_* R\cHom_Y(\cG,Rj^!L).\]
Il suffit d'appliquer \eqref{eq.Rf!1} et \eqref{eq.Rf*1} (\resp
\eqref{eq.Rf!3} et \eqref{eq.Rf*3}).
\end{proof}
\end{ssect}

\section{Variantes et cycles proches}\label{S5}
\begin{vari}\label{ss.vari.fini}
Soient $k$ un corps fini, $l$ un nombre premier $\neq \car(k)$. Soit
$X$ un \sch \tf sur $k$. Soit $r\in \Q$. On fixe un plongement
$\iota : \Qb \to \Qlb$. On dit qu'un $\Qlb$-\faisc $\cF$ est
$r$-\emph{entier} (\resp $r$-\emph{\entinv{}}) si pour tout point
\geom $\xb$ au-dessus d'un point fermé $x$ de $X$, et toute valeur
propre $\al$ du Frobenius \geom $F_x \in
\Gal(\kappa(\xb)/\kappa(x))$ agissant sur $\cF_\xb$,
$\al/\iota(q^r)$ (\resp $\iota(q^r)/\al$) est entier sur $\Z$, où
$q=\sharp \kappa(x)$. Cette définition ne dépend pas du choix de
$\iota$. On définit l'intégralité pour $K\in \Dbc(X,\Qlb)$ de
manière analogue à \ref{defn.c}.

On a des résultats similaires pour les diverses opérations :
\ref{theo.Rf!b} à \ref{theo.RHom}, \ref{prop.Rj*dcn},
\ref{prop.Rf*qf}, \ref{prop.i^!}. Le cas entier de l'analogue de
\ref{theo.Rf!b} est un théorème de Deligne \cite[XXI 5.2.2]{SGA7}.
Les démonstrations des autres résultats sont similaires à celles
données aux \S\S\ 3 et 4.
\end{vari}

\begin{vari}\label{vari.Xse}
Soient $R$ un anneau de valuation discrète hensélien excellent de
corps résiduel fini $k = \F_{p^\nu}$, $K$ son corps des fractions,
$S=\Spec R$, $\eta= \Spec K$, $s=\Spec k$.  Soit $X$ un \sch \tf sur
$s$. On a un topos $X\times_s \eta$ \cite[XIII 1.2.4]{SGA7}.
Rappelons qu'un faisceau d'ensembles sur $X\times_s \eta$ est un
faisceau sur $X_\sbar$ muni d'une action continue \ibid[1.1.2] de
$\Gal(\Kb/K)$, compatible à l'action de $\Gal(\Kb/K)$ sur $X_\sbar$
(via $\Gal(\Kb/K) \to \Gal(\kb/k)$). Fixons un nombre premier $l
\neq p$. Soit $\cF\in \Mod_c(X\times_s \eta,\Qlb)$. Pour $x\in X$,
soit $\Phi_x\in \Gal(\kb/\kappa(x)) \times_{\Gal(\kb/k)}
\Gal(\Kb/K)$ un relèvement du Frobenius géométrique $F_x\in
\Gal(\kb/\kappa(x))$. D'après le théorème de monodromie locale,
les valeurs propres de $\Phi_x$ agissant sur $\cF_\xb$ sont bien
définies à multiplication près par des racines de l'unité. Soit
$r\in \Q$. On fixe $\iota : \Qb \to \Qlb$. On dit que $\cF$ est
\emph{$r$-entier} (\resp \emph{$r$-entier inverse}) si pour tout
$x\in |X|$ et toute valeur propre $\al$ de $\Phi_x$ agissant sur
$\cF_\xb$, $\al/\iota(q^r)$ (\resp $\iota(q^r)/\al$) est entier sur
$\Z$, où $q=\sharp \kappa(x)$. Cette définition ne dépend pas des
choix de $\Phi_x$ et de $\iota$. On définit l'intégralité pour
$K\in \Dbc(X\times_s\eta,\Qlb)$ de manière analogue à
\ref{defn.c}.

Toute section continue $\sigma$ de $\Gal(\Kb/K) \to \Gal(\kb/k)$
induit un foncteur exact
\[ \sigma^* : \Dbc(X\times_s \eta,\Qlb) \to \Dbc(X,\Qlb). \]
Un complexe $K\in \Dbc(X\times_s \eta,\Qlb)$ est \eent (\resp
\eentinv) si et seulement si $\sigma^* K$ l'est. Comme $\sigma^*$
commute aux six opérations et à la dualité, on déduit de
\ref{ss.vari.fini} des résultats similaires pour ces opérations :
\ref{theo.Rf!b} à \ref{theo.RHom}, \ref{prop.Rj*dcn},
\ref{prop.Rf*qf}, \ref{prop.i^!}.
\end{vari}

Soient $S$, $\eta$, $s$, $l$ comme dans \ref{vari.Xse}. Le résultat
principal de ce~\S\  est le suivant.

\begin{theo}\label{theo.RP}
Soit $X$ un \sch \tf sur $S$. Le foncteur des cycles proches
\[R\Psi_X : \Dbc(X_\eta,\Qlb) \to \Dbc(X_s\times_s \eta,\Qlb)\]
induit
\begin{align*}
\Dbclent{X_\eta} &\to \Dbclent{X_s\times_s \eta},\\
\Dbcleentinv[1,0]{X_\eta} &\to \Dbcleentinv[1,0]{X_s\times_s \eta}.\\
\end{align*}
\end{theo}

\begin{ssect}\label{ss.S}
Soient $S=\Spec R$ un trait hensélien quelconque, $\eta$ son point
générique, $s$ son point fermé. On va garder ces notations jusqu'en
\ref{lemm.RP}.
\end{ssect}

\begin{defn}
(a) Soient $X$ un $S$-\sch \tf, $Z$ une partie fermée contenant
$X_s$. On dit que le couple $(X,Z)$ est \emph{\sst{}} si, localement
pour la \topet, il est de la forme
\[(\Spec R[t_1, \dots ,t_n]/(t_1 \dots t_r-\pi),Z),\] où $\pi$ est une uniformisante de $R$,
$Z$ est défini par l'idéal $(t_1 \dots t_s)$, $1\le r\le s \le n$.
Le couple est dit \emph{\sss{}} s'il est \sst et si $Z$ est la somme
d'une famille finie de diviseurs \regs de $X$.

(b) Soit $X$ un $S$-\sch \tf. On dit que $X$ est \emph{\sss{}} si
$(X,X_s)$ est un couple \sss.
\end{defn}

Soit $(X,Z)$ un couple \sss avec $Z=\sum_{i\in I} D_i$, $(D_i)_\idI$
une famille finie de diviseurs \regs. Alors $X_s=\bigcup_{i\in I-J}
D_i$, où $J= \ensdr{i\in I}{D_i \not\subset X_s}$. Soit $H=
\bigcup_{j\in J} D_j$ la réunion des composantes horizontales.
Alors $Z=X_s\cup H$.

Nous établirons d'abord \ref{theo.RP} dans le cas \sst. Nous aurons
besoin pour cela des points \tmii et \tmiii du lemme suivant (la
partie \tmi est utilisée dans la preuve de (ii)) :

\begin{lemm}\label{lemm.RPdev}
Soient $(X,Z)$ un couple \sss 
sur $S$, $Z=\sum_{i\in I} D_i$ avec $(D_i)_\idI$ une famille finie
de diviseurs \regs, $J$ comme plus haut, $U=X-Z$, $u$ l'inclusion
$U\hra X_\eta$, $\Lbd = \Z/m \Z$ avec $m$ inversible sur $S$, $\cG
\in \Mod_c(U,\Lbd)$ \loccst, \mr sur~$X$.

\tmi Soient $i\in J$, $U_\ip = X-\bigcup_{h\in I-\{i\}} D_h$, $D_{i,
U_\ip} = D_i \times_X U_\ip$, d'où un diagramme à carrés
cartésiens
\[\xymatrix{&D_{i,U_\ip}\ar@{^{(}->}[r]^{j'_{\ip}}\ar[d]^{\iota'_{i}}
& (D_i)_\eta \ar@{^{(}->}[r] \ar[d]^{(\ii)_\eta} &
D_i\ar[d]^{\ii} & (D_i)_s \ar[l] \ar[d]^{(\ii)_s}\\
U\ar@{^{(}->}[r]^{j^\ip} \ar@/_1pc/[rr]_{u} &
U_\ip\ar@{^{(}->}[r]^{j_\ip}& X_\eta \ar@{^{(}->}[r] & X&
X_s\ar[l]}\] Alors la flèche
\beq
\al : (\ii)_s^* R\Psi_X Rj_{\ip*} (Rj^\ip_* \cG)  \to R\Psi_{D_i}
Rj'_{\ip*} {\ii'}^{*} (Rj^\ip_* \cG)
\eeq
composée de $(\ii)_s^* R\Psi_X (Ru_* \cG) \to R\Psi_{D_i}
(\ii)_\eta^* (Ru_* \cG)$ \cite[XIII (2.1.7.2)]{SGA7} et du \chgt de
base $R\Psi_{D_i} (\ii)_\eta^* Rj_{\ip *} (Rj^\ip_* \cG) \to
R\Psi_{D_i} Rj'_{\ip*} {\ii'}^{*} (Rj^\ip_* \cG)$  est un \isom.

\tmii 
Soient $i\in I-J$, $U_\ip$, $D_{i,U_\ip}$ comme plus haut, d'où un
diagramme à carrés cartésiens
\[\xymatrix{U\ar@{^{(}->}[r] \ar@{^{(}->}[dd]^{u} &U_\ip \ar@{^{(}->}[dd]
& D_{i,U_\ip}\ar[l]\ar@{^{(}->}[d]^{j'_\ip}\\
&& D_i \ar[d]^{(\ii)_s}\\
X_\eta \ar@{^{(}->}[r] & X & X_s \ar[l]}\] Alors $R\Psi_{U_\ip} \cG
\simeq \Psi_{U_\ip} \cG \in \Mod_c(D_{i,U_\ip}\times_s \eta,\Lbd)$
est lisse, \mr sur $D_i$, et le morphisme
\beq\label{eq.RPii}
\beta : (\ii)_s^* R\Psi_X Ru_*\cG \to Rj'_{\ip*} R\Psi_{U_\ip} \cG
\eeq
déduit de $R\Psi_X Ru_* \cG \to (\ii)_{s*} Rj'_{\ip*} R\Psi_{U_\ip}
\cG$ \cite[XIII (2.1.7.1)]{SGA7} est un \isom.

\tmiii Supposons que $J=\emptyset$
. Soit $f:Y\to X$ un morphisme de \schs \tq $Y$ soit un $S$-\sch
\sss avec $(f^{-1}(D_i))_\idI$ une famille de diviseurs \regs
de~$Y$. Alors le morphisme
\beq\label{eq.RPiii}
\gamma : f_s^* R\Psi_X \cG \to R\Psi_Y f_\eta^* \cG
\eeq
\cite[XIII
(2.1.7.2)]{SGA7} est un \isom.
\end{lemm}

Le point (ii) est une généralisation partielle de \cite[1.5
(a)]{IllPL}.

\begin{proof}
On remplace tout d'abord le premier énoncé de (ii) par l'assertion
que les faisceaux $R^q \Psi_{U_{(i)}} \cG$, $q\in \Z$ sont lisses,
modérément ramifiés sur $D_i$. L'annulation des cycles proches
supérieures résultera de (iii).

La question est locale sur $X$. Soit $y$ un point de $X_s$. \Ops que
$y\in D_j$ pour tout $j\in I$. Soit $\pi$ une uniformisante de $R$.
Il existe un ouvert de~$X$ contenant~$y$, lisse sur
\[\Spec
R[t_j]_{j\in I} / (\prod_{j\in I - J} t_j - \pi)
\]
avec $t_j$ définissant $D_j$. En vertu du lemme d'Abhyankar, il
existe, au voisinage de~$y$, un revêtement fini $g: \Xt =
X[T_j]_{j\in I}/(T_j^{n}-t_j)_{j\in I} \to X$ où $n$ est un entier
premier à l'exposant \cara de $s$, \tq $(g|U)^*\cG$ se prolonge en
un $\Lbd$-module \loct \cst \constr sur~$\Xt$. Comme $\cG$ s'injecte
dans $(g|U)_* (g|U)^* \cG$ et le quotient $\cG_1$ est \mr sur $X$,
on peut itérer cette construction. Pour tout $N\ge 1$, on obtient,
quitte à rétrécir $X$, une résolution
\[\cG \to  (g|U)_* (g|U)^*
\cG \to (g_1|U)_* (g_1|U)^* \cG_1 \to \dots \to (g_N|U)_* (g_N|U)^*
\cG_N.
\]
Donc on est ramené à montrer le lemme pour le faisceau $(g|U)_*
(g|U)^* \cG$.

{\sloppy Soient $R_1$ l'anneau obtenu en adjoignant à $R$ les
$n$-ièmes racines de l'unité, $R'= R_1[\Pi]/{(\Pi^n-\pi)}$,
$S'=\Spec R'$. Comme $R\Psi$ commute au \chgt de traits $S'\to S$
\cite[Th.\ finitude, 3.7]{SGA4d}, il suffit de montrer le lemme pour
$((g|U)_* (g|U)^* \cG)_{S'} \simeq (g_{S'}|U_{S'})_{*}
(g_{S'}|U_{S'})^* \cG_{S'}$. Or $X_{S'} = \coprod_\zeta X_\zeta$,
où $X_\zeta$ est lisse sur
\[\Spec R[T_j]_{j\in J} / (\prod_{j\in I - J} T_j - \zeta \cdot \Pi),\]
$\zeta$ parcourt les $n$-ièmes racines de l'unité. Donc $(X_{S'},
g_{S'}^{-1}(Z_{S'})_\red)$ est un couple \sss avec
$(g_{S'}^{-1}(D_{i,S'})_\red)_{i\in I}$ une famille de diviseurs
\regs.

}On peut alors se ramener à montrer le lemme pour un \faisc $\cG$
qui se prolonge en un $\Lbd$-module \loct constant \constr sur $X$,
puis au cas $\cG = \Lbd_U$ par la formule de projection.

L'assertion (iii) découle alors de la fonctorialité de
\cite[3.3]{Illss}. Plus précisément, on a les diagrammes \comms
\[\xymatrix{\wedge^q f_s^* R^1\Psi_X \Lbd_{X_\eta} \ar[r]^{\sim}\ar[d]^{\wedge^q \fH^1 \gamma}
& f_s^* R^q\Psi_X \Lbd_{X_\eta}\ar[d]^{\fH^q \gamma} \\
\wedge^q R^1\Psi_Y f_\eta^* \Lbd_{X_\eta} \ar[r]^{\sim} & R^q\Psi_Y
f_\eta^*\Lbd_{X_\eta}}\] et
\[\xymatrix{0\ar[r]& f_s^*\Lbd_{X_s}(-1) \ar[r]\ar[d]^{\simeq}
& f_s^* \iota_X^* R^1 j_{X*}\Lbd_{X_\eta} \ar[r]\ar[d]^{\ref{lemm.cb} \tmii}_\simeq & f_s^*
R^1\Psi_X\Lbd_{X_\eta}\ar[d]^{\fH^1{\gamma}} \ar[r]& 0\\
0\ar[r] & \Lbd_{Y_s}(-1) \ar[r] & \iota_Y^* R^1 j_{Y*} f_\eta^*
\Lbd_{X_\eta}\ar[r] & R^1\Psi_Y f_\eta^* \Lbd_{X_\eta} \ar[r]& 0}\]
où $j_X:X_\eta \hra X$, $j_Y: Y_\eta \hra Y$, $\iota_X : X_s \to
X$, $\iota_Y : Y_s \to Y$. Les lignes du deuxième diagramme sont
des suites exactes courtes et le carré à gauche est donné par le
\diacom
\[\xymatrix{f_s^*\Lbd_{X_s}(-1) \ar[r]^-d\ar[d]^{\simeq}
& f_s^* \bigoplus_{\idI}\Lbd_{X_i}(-1)\ar[r]_c^{\sim} \ar[d]^{\simeq}
&f_s^* i_X^* R^1 j_{X*}\Lbd_{X_\eta} \ar[d]\\
\Lbd_{Y_s}(-1) \ar[r]^-d &
\bigoplus_{\idI}\Lbd_{Y_i}(-1)\ar[r]^{\sim}_c&i_Y^* R^1 j_{Y*}
f_\eta^* \Lbd_{X_\eta}}\] où les flèches marquées $d$ sont des
diagonales et celles marquées $c$ sont induites par les classes des
diviseurs \regs.

(i) Il s'agit de montrer l'assertion suivante :
\begin{cond}
(A) Le morphisme de foncteurs
\beq\label{eq.RPfonct}
(\ii)_s^* R\Psi_X Rj_{\ip*} \to R\Psi_{D_i} Rj'_{\ip*} {\ii'}^{*}
\eeq
induit un \isom sur $Rj^\ip_* \Lbd_U$.
\end{cond}

On a un triangle \dist
\[(\ii')_* \Lbd_\DiUi (-1)[-2] \to \Lbd_{U_\ip} \to Rj^\ip_* \Lbd_U \to. \]
Comme \eqref{eq.RPfonct} induit trivialement un \isom sur le premier
terme, (A) équivaut à
\begin{cond}
(B) Le morphisme $\text{\eqref{eq.RPfonct}} \Lbd_{U_\ip}$
\beq\label{eq.RPmod}
(\ii)_s^* R\Psi_X Rj_{\ip*}  \Lbd_{U_\ip} \to R\Psi_{D_i} Rj'_{\ip*}
\Lbd_\DiUi
\eeq
est un \isom.
\end{cond}

On montre ces énoncés par \rec sur $\sharp J \ge 1$. Le cas
$\sharp J = 0$ est vide.

Dans le cas \gel, on montre d'abord que pour tout $j\in J-\{i\}$,
$\text{\eqref{eq.RPmod}} | (D_{ij})_s$ est un \isom, où $D_{ij} =
D_i \cap D_j$. Soit $U_\ijp = X-\bigcup_{h\in I-\{i,j\}} D_h$. On
considère le diagramme commutatif
\[\xymatrix{
\DiUi \ar@{^{(}->}[rr]^{j'_1}\ar[dd]^{\ii'} & & D_{i,U_\ijp}
\ar@{-}[d] \ar@{^{(}->}[rr]^{j_{2,i}} && (D_i)_\eta \ar@{-}[d]
\ar@{^{(}->}[rr] && D_i\ar[dd]^{\iota_{i}}\\
& D_{ij,U_\ijp} \ar@{->}[ru]\ar[dd]\ar@{^{(}->}[rr]^(.3){j'_2} &
\ar@{-->}[d] & (D_{ij})_\eta \ar@{^{(}->}[rr]\ar[dd]\ar@{->}[ru] &
\ar@{-->}[d] & D_{ij}
\ar[dd]_(.7){\iota_{i,j}}\ar[ru]^{\iota_{j,i}}\ar[rd]^{\iota_{ij}} & \\
U_\ip \ar@{^{(}-}[r]^(.6){j_1} &\ar@{-->}[r]& U_\ijp
\ar@{^{(}-->}[rr]^(.3){j_2} && X_\eta
\ar@{^{(}-->}[rr] && X\\
& D_{j, U_\ijp} \ar@{-->}[ru] \ar@{^{(}->}[rr]^{j_{2,j}} &&
(D_j)_\eta \ar@{-->}[ru] \ar@{^{(}->}[rr] && D_j \ar[ru]_{\iota_j}
&}
\]
D'après \cite[XII 4.4 (i)]{SGA4}, le composé
\begin{multline*}
(\iota_{ij})_s^* R\Psi_X R(j_2j_1)_* \Lbd_{U_\ip}\\
\begin{aligned}
\to & (\iota_{j,i})_s^* R\Psi_{D_i} R(j_{2,i} j'_{1})_* \Lbd_{\DiUi}
&& \text{\eqref{eq.RPmod}}
| (D_{ij})_s\\
\simto & R\Psi_{D_{ij}} Rj'_{2*} ((Rj'_{1*} \Lbd_{\DiUi})|
D_{ij,U_\ijp}) && \text{\hyp de \rec (A)}
\end{aligned}
\end{multline*}
est égal au composé
\begin{multline*}
(\iota_{ij})_s^* R\Psi_X R(j_2j_1)_* \Lbd_{U_\ip}\\
\begin{aligned}
\simto & (\iota_{i,j})_s^* R\Psi_{D_j} R(j_{2,j})_*
((Rj_{1*}\Lbd_{U_\ip})|D_{j,U_\ijp}) &&
\text{\hyp de \rec (A)}\\
\to & R\Psi_{D_{ij}}Rj'_{2*}((Rj_{1*} \Lbd_{U_\ip})|D_{ij,U_\ijp}) &&\text{($*$)}\\
\simto & R\Psi_{D_{ij}}Rj'_{2*}((Rj'_{1*}
\Lbd_{\DiUi})|D_{ij,U_\ijp}) && \text{\ref{lemm.cb} \tmii}
\end{aligned}
\end{multline*}
où ($*$) est un morphisme de type \eqref{eq.RPfonct} appliqué à
$(Rj_{1*} \Lbd_{U_\ip})| D_{j,U_\ijp}$. On a le triangle \dist
\[\Lbd_{D_{j,U_\ijp}}(-1)[-2] \to \Lbd_{D_{j,U_\ijp}} \to (Rj_{1*} \Lbd_{U_\ip})| D_{j,U_\ijp}
\to.
\]
Donc l'\hyp de \rec (B) implique que ($*$) est un \isom. Il en
résulte que $\text{\eqref{eq.RPmod}}| (D_{ij})_s$ est un \isom.

Il reste à montrer que $\text{\eqref{eq.RPmod}}|(V_i)_s$ est un
\isom, où
\[V_i = X-\bigcup_{h\in J-\{i\}} D_h.\] Comme $(V_i)_\eta = U_\ip$, $j_{\ip, V_i} = \id_{U_\ip}$, ceci découle de (iii).

(ii) Comme $U_\ip$ est lisse sur $S$, $R\Psi_{U_\ip} \Lbd_U \simeq
\Lbd_{\DiUi}$, donc est \mr sur $D_i$.

Pour montrer que $\beta$ est un \isom, traitons d'abord deux cas
spéciaux : (a) $\sharp J = 0$ ; (b) $\sharp (I-J) = \sharp J = 1$.

Dans le cas (a), on a $U=X_\eta$, $j_\eta=\id_{X_\eta}$. On pose
$D=D_i$, $E=\bigcup_{j\in I-\{i\}}D_j$, $D^*=D-D\cap E=\DiUi$, d'où
un \diacom
\[\xymatrix{U\ar@{^{(}->}[r]^{j^\ip} \ar@{=}[dd] & X-E \ar@{^{(}->}[dd]^{j_\ip}
& D^* \ar[l]_{\ii'} \ar@{^{(}->}[d]^{j'_\ip}\\
&& D\ar[ld]_{\ii} \ar[d]^{(\ii)_s}\\
X_\eta \ar@{^{(}->}[r]^j & X & X_s \ar[l]}
\]
On a le carré commutatif
\[\xymatrix{\wedge^q(\ii)_s^* R^1\Psi_X \Lbd_U \ar[r]^{\sim}\ar[d]^{\wedge^q \fH^1 \beta} &
(\ii)_s^* R^q\Psi_X \Lbd_U \ar[d]^{\fH ^q \beta}\\
\wedge^q R^1 j'_{\ip*} R\Psi_{U_\ip} \Lbd_U \ar[r] &
R^qj'_{\ip*} R\Psi_{U_\ip} \Lbd_U\\
\wedge^q R^1 j'_{\ip*} \Lbd_{D^*} \ar[u]^{\simeq} \ar[r]^{\sim} &
R^qj'_{\ip
*} \Lbd_{D^*} \ar[u]^{\simeq}}
\]
Donc il suffit de montrer que $\fH^1\beta$ est un \isom.

On a le \diacom
\[\xymatrix{& \ii^* Rj_* \Lbd_U \ar@{=}[d] \ar[r]^{p_1} &
(\ii)_s^* R\Psi_X \Lbd_U \ar[dd]^{\beta} \\
\ii^* Rj_{\ip *}\Lbd_{X-E} \ar[d]^{\simeq}_{b_1} \ar[r]^{r_1} &
\ii^* Rj_{\ip*}Rj^\ip_*\ar[d]_{b_2}
\Lbd_U &\\
Rj'_{\ip*} \ii^{'*} \Lbd_{X-E} \ar[r]^{r_2} \ar@/_2pc/[rr]^{\sim} &
Rj'_{\ip*} \ii^{'*}Rj^\ip_* \Lbd_U \ar[r]^{p_2} & Rj'_{\ip*}
R\Psi_{X-E} \Lbd_U}
\]
où $b_1$, $b_2$ sont des changements de base, $r_1$, $r_2$ sont
induits par l'adjonction $\Lbd_{X-E} \to Rj^\ip_* \Lbd_U$. La flèche
$b_1$ est un isomorphisme en vertu de \ref{lemm.cb} (ii). Le composé
$p_2 r_2$ est induit de l'\isom $\Lbd_{D^*} \simto R\Psi_{X-E}
\Lbd_U$, donc est un \isom. On a
\begin{align*}
\fH^1 \ii^* Rj_{\ip *}\Lbd_{X-E} &= \bigoplus_{j\in
I-\{i\}}(\iota_{j,i})_* \Lbd_{D_{ij}},\\
\fH^1 \ii^* Rj_* \Lbd_U &= \Lbd_D \oplus \bigoplus_{j\in I-\{i\}}
(\iota_{j,i})_* \Lbd_{D_{ij}},
\end{align*}
$\fH^1 (\ii)_s^* R\Psi_X \Lbd_U$ est le quotient de $\Lbd_D \oplus
\bigoplus_{j\in I-\{i\}} (\iota_{j,i})_* \Lbd_{D_{ij}}$ par
$\Lbd_{D_i}$ inclus diagonalement, $\fH^1 r_1$ est l'inclusion dans
le second membre, $\fH^1 p_1$ est la projection. Donc $\fH^1(p_1
r_1)$ est un \isom. Il s'en suit que $\fH^1 \beta$ est un \isom.
D'où (a).

Dans le cas (b), on a $D_i = X_s$, $(\ii)_s = \id_{X_s}$, $H=D_j$ où
$j$ est l'\elt de $J$. On pose $V= U_\ip = X-H$, d'où un diagramme à
carrés cartésiens
\[\xymatrix{& U \ar@{^{(}->}[d]^{u}\ar@{^{(}->}[rr] && V\ar@{^{(}->}[d]^v &&
V_s
\ar@{^{(}->}[d]^{v_s}\ar[ll]_{i_V} \\
& X_\eta\ar@{^{(}->}[rr] && X && X_s \ar[ll]_{i}\\
H_\eta\ar[ru] \ar@{^{(}->}[rr] && H \ar[ru]^{h} && H_s
\ar[ru]\ar[ll]}
\]
On a le triangle \dist
\[i^* Rv_* \Lbd_V \sto{\beta'} R\Psi_X Ru_* \Lbd_U \to R\Phi_X Rv_* \Lbd_V \to.\]
On a
\[R^q v_* \Lbd_V =
\begin{cases}
\Lbd_X & \text{si $q=0$,}\\
h_* \Lbd_H (-1) & \text{si $q=1$,}\\
0 & \text{sinon.}
\end{cases}
\]
Donc $R\Phi_X R^q v_* \Lbd_V = 0$, pour tout $q$. Donc $R\Phi_X Rv_*
\Lbd_V = 0$, $\beta'$ est un \isom. Par ailleurs, on a le \diacom
\[\xymatrix{i^* Rv_* \Lbd_V \ar[r]^{\beta'} \ar[d]^{\simeq}_{\text{\ref{lemm.cb} (ii)}}
& R\Psi_X Ru_* \Lbd_U\ar[d]^{\beta} \\
Rv_{s*} i_V^* \Lbd_V \ar[r]^{\sim} & Rv_{s*}R\Psi_V \Lbd_U}
\]
Donc $\beta$ est un \isom. D'où (b).

Pour le cas général, procédons par \rec sur $\sharp J$. Le cas
$\sharp J=0$ est le cas (a) traité plus haut. Supposons $\sharp J
\ge 1$. Prenons $j\in J$, d'où un \diacom
\[\xymatrix{&U\ar@{^{(}->}[rr]\ar@{^{(}->}[d]^{j^\jp} && U_\ip\ar@{^{(}->}[d]
&& \DiUi\ar[ll]\ar@{^{(}->}[d]^{j_1}\\
&U_\jp\ar@{^{(}-}[d]\ar@{^{(}->}[rr] && U_\ijp\ar@{^{(}-}[d] &&
D_{i,U_\ijp}\ar[ll]\ar@{^{(}->}[dd]^{j_2}\\
\DjUj\ar@{^{(}->}[dddd]^{u'_\jp} \ar@{^{(}->}[rr]\ar[ur]^{\iota'_j}
&\ar@{-->}[ddd]& D_{j,U_\ijp}\ar@{^{(}->}[dddd]\ar[ur]
&\ar@{-->}[ddd]&
D_{ij,U_\ijp}\ar[ll]\ar[ur]^{\iota'_{j,i}}\ar@{^{(}->}[dd]^{j'_2} & \\
&&&&&D_i\ar[dd]^{(\iota_i)_s}\\
&&&&D_{ij}\ar[dd]_(0.7){(\iota_{i,j})_s}\ar[ur]^{\iota_{j,i}}\\
&X_\eta \ar@{^{(}-->}[rr] && X && X_s\ar@{-->}[ll]\\
(D_j)_\eta\ar@{-->}[ur]^{(\iota_j)_\eta} \ar@{^{(}->}[rr] &&
D_j\ar@{-->}[ur]^{\iota_j} && (D_j)_s \ar[ll]\ar[ur]_{(\iota_j)_s}
}
\] On a un \diacom dans $\Dbc(D_{ij}\times_s\eta,\Lbd)$
\beq\tag{$\dagger$}
\xymatrix@C=0.5cm{\iota_{j,i}^* (\iota_i)_s^* R\Psi_X Ru_* \Lbd_U
\ar[r]^{\beta_1} \ar@{=}[d] & \iota_{j,i}^* Rj_{2*}
R\Psi_{U_\ijp}Rj^\jp_* \Lbd_U
\ar[r]^{\text{(b)}}_\sim\ar[d]^{\text{chg\up{t} de base}} &
\iota_{j,i}^* R{j_2j_1}_* R\Psi_{U_\ip} \Lbd_U
\ar[d]^{\text{\ref{lemm.cb} (i)}}_\simeq\\
(\iota_{i,j})_s^* (\iota_j)_s^* R\Psi_X Ru_* \Lbd_U
\ar[d]^{\tmi}_\simeq & Rj'_{2*} {\iota'_{j,i}}^*
R\Psi_{U_\ijp}Rj^\jp_* \Lbd_U
\ar[r]^{\text{(b)}}_\sim\ar[d]^{\tmi}_\simeq & Rj'_{2*}
{\iota'_{j,i}}^*
Rj_{1*} R\Psi_{U_\ip} \Lbd_U\\
(\iota_{i,j})_s^* R\Psi_{D_j} Ru'_{\jp*} {\iota'_j}^* Rj^\jp_*
\Lbd_U \ar[r]^-{\beta_2} & Rj'_{2*} R\Psi_{D_{j,U_\ijp}}{\iota'_j}^*
Rj_* \Lbd_U}
\eeq
où $\beta|D_{ij}$ est le composé des deux flèches de la première
ligne de ($\dagger$), $\beta_1$ est induit par une flèche de type
\eqref{eq.RPii} et $\beta_2$ est une flèche de type \eqref{eq.RPii}.
La commutativité du carré à droite est claire et celle du carré à
gauche se voit en appliquant \ref{lemm.RP} au carré
\[\xymatrix{
D_{j,U_\ijp} \ar[r]\ar@{^{(}->}[d] & U_\ijp \ar@{^{(}->}[d] \\
D_j \ar[r] & X}
\]
La flèche $\beta_2$ est un \iso en vertu de l'\hyp de \rec et du
\trdist
\[\Lbd_{\DjUj}(-1)[-2] \to \Lbd_{\DjUj} \to {\iota'_j}^* Rj_* \Lbd_U \to.\]
Donc $\beta|D_{ij}$ est un \iso. Il reste à montrer que
$\beta|(D_i-D_{ij})$ est un \iso, ce qui résulte de l'\hyp de \rec.
\end{proof}

\begin{lemm}\label{lemm.RP}
Soient
\[\xymatrix{
X'\ar[r]^h \ar[d]^{f'} & X\ar[d]^f \\
Y'\ar[r]^g & Y}
\]
un carré \comm de $S$-\schs, $\Lbd$ un anneau. Alors on a un \diacom
de foncteurs $D^+(X_\eta,\Lbd) \to D^+(Y'_s \times_s \eta,\Lbd)$
\[\xymatrix{
& g_s^* Rf_{s*} R\Psi_X \ar[r] & Rf'_{s*} h_s^* R\Psi_X \ar[dr] \\
g_s^* R\Psi_Y Rf_{\eta*} \ar[ur] \ar[dr] &&& Rf'_{s*} R\Psi_{X'}
h_\eta^* \\
& R\Psi_{Y'} g_\eta^* Rf_{\eta *} \ar[r] & R\Psi_{Y'} Rf'_{\eta *}
h_\eta^* \ar[ur]}
\]
où les flèches horizontales sont des changements de base, les
flèches montantes sont \cite[XIII (2.1.7.1)]{SGA7}, les flèches
descendantes sont \ibid[XIII (2.1.7.2)].
\end{lemm}

\begin{proof}
Soient $K= \kappa(\eta)$, $\Kb$ une clôture \sepble de $K$, $\etab =
\Spec \Kb$, $\Sbar$ le normalisé de $S$ dans $\etab$. On ajoute une
barre au-dessus pour le changement de base $\lbar{S} \to S$. On note
$i: \sbar \to \Sbar$, $j: \etab\to \Sbar$. Il suffit de montrer la
commutativité du \diag de foncteurs $D^+(X_{\etab}, \Lbd) \to
D^+(Y'_{\sbar}, \Lbd)$
\catcode`!12
\[\xymatrix@!C0@C=4.5em@R=8ex{
&& g_{\sbar}^* Rf_{\sbar *} i_X^* Rj_{X*} \ar[rr] && Rf'_{\sbar*}
h_{\sbar}^* i_X^* Rj_{X*} \ar@{-}[rd]^{\sim} \\
& g_{\sbar}^* i_Y^* R\bar{f}_* Rj_{X*} \ar[ur] \ar@{-}[dr]^{\sim}
&&&& Rf'_{\sbar*}i_{X'}^*
\bar{h}^* Rj_{X*} \ar[rd] \\
g_{\sbar}^* i_Y^* Rj_{Y*} Rf_{\etab *} \ar@{-}[ur]^{\sim}
\ar@{-}[dr]^{\sim} && i_Y^* \bar{g}^* R\bar{f}_* Rj_{X*} \ar[rr] &&
i_Y^* R\lbar{f'}_* \bar{h}^* Rj_{X*} \ar[ur] \ar[dr] && Rf'_{\sbar*}
i_{X'}^* Rj_{X'*}
h_{\etab}^* \\
& i_{Y'}^* \bar{g}^* Rj_{Y*} Rf_{\etab*} \ar@{-}[ur]^{\sim} \ar[dr]
&&&& i_{Y'}^* R\lbar{f'}_* Rj_{X'*} h_{\etab}^* \ar[ur] \\
&& i_{Y'}^* Rj_{Y'*} g_{\etab}^* Rf_{\etab*} \ar[rr] && i_{Y'}^*
Rj_{Y'*} Rf'_{\etab*} h_{\etab}^* \ar@{-}[ur]^{\sim} }
\]\catcode`!\active
où toutes les flèches sont des changements de base. La commutativité
de la cellule en haut (\resp en bas) résulte de \cite[XII 4.4
(i)]{SGA4} (\resp \cite[XII 4.4 (ii)]{SGA4}). Les commutativités des
deux autres cellules sont triviales.
\end{proof}

\begin{prop}\label{prop.RPss}
Soient $S$ comme dans \ref{vari.Xse}, $(X,Z)$ un couple \sst
sur~$S$, $U= X-Z$, $u: U\to X_\eta$, $\cG \in \Mod_c(U,\Qlb)$ \ent
(\resp \entinv) lisse, \mr sur~$X$. Alors $R\Psi_X Ru_* \cG$ est
\eent[1,0] (\resp \eentinv[1,0]).
\end{prop}

\begin{proof} \Ops que $(X,Z)$ est un couple \sss sur~$S$.

Traitons d'abord le cas particulier où $Z$ est un diviseur \reg.
Alors $Z=X_s$, $U=X_\eta$, $u=\id$, $X$ est lisse sur $S$. Soit
$x\in |X_s|$. Quitte à faire un changement de traits étale, \ops que
$X\to S$ admet une section $\sigma$ \tq $\sigma(s)=x$. D'après
\ref{lemm.RPdev} (iii), $(R\Psi_X \cG)_x = \sigma_s^*R\Psi_X\cG
\simeq R\Psi_S \sigma_\eta^*\cG$ est $I$-entier (\resp $I$-entier
inverse), car $R\Psi_S$ s'identifie à l'identité.

Le cas \gel découle de \ref{lemm.RPdev}~(ii), du cas spécial
ci-dessus, et de la variante \ref{vari.Xse} de \ref{prop.Rj*dcn}
(\resp \ref{prop.Rj*dcninv}) au-dessus de $s\times_s \eta$.
\end{proof}


\ssect

Pour la démonstration de \ref{theo.RP}, nous aurons besoin du lemme
\ref{lemm.superdJ} ci-après, analogue de \ref{lemm.dJ}.

Soient $S$ un trait hensélien excellent, $X$ un $S$-\sch \sep \tf,
$U\subset X_\eta$ une partie ouverte. Pour $f : S'\to S$ un
morphisme fini de traits et $g : Y\to X_{S'}$ un morphisme propre de
\schs, on considère la condition suivante :
\begin{numcond}\label{cond.superdJ}
On a $Y= Y_1 \coprod Y_2$, où $Y_1$ est strictement semi-stable
sur~$S'$, $g^{-1}(U) \subset Y_2$ et $(Y_2,Y_2 - g^{-1}(U))$ est un
couple strictement semi-stable sur~$S'$.
\end{numcond}

\begin{lemm}\label{lemm.superdJ}
\tmi Il existe un morphisme fini de traits $S'\to S$ et un morphisme
propre $r_0 : X_0' \to X_{S'}$ de \schs \ver ant \ref{cond.superdJ}
(où $Y= X_0'$) et \tsq $(r_0)_\eta$ surjectif.

\tmii Pour $n \ge 0$, il existe un morphisme fini de traits $f :
S'\to S$ et une augmentation de \sch simplicial $n$-tronqué \sscin
$r\sm : X'\sm \to X_{S'}$ \tsq pour $0\le m\le n$, $f$ et $r_m$ \ver
ent \ref{cond.superdJ} (où $Y= X_m'$) et que ${r\sm}_\eta$ soit un
\hyperr propre $n$-tronqué.
\end{lemm}

\begin{proof}
(i) Cas $X$ intègre et $X_\eta$ \geot \irr. Résulte de
\cite[6.5]{deJong}. \vcb{v1m}{Notons que l'\hyp dans \ibid que $S$
soit complet peut être remplacée par l'excellence de~$S$, voir
\cite[3.8]{Zheng}.}

Cas \gel. \Ops que les composantes \irrs de $X_\eta$ sont \geot
\irrs. On fait une \rec sur le nombre $n$ de composantes \irrs de
$X_\eta$.

Si $n=0$, alors $X_\eta$ est vide. On prend $S'=S$, $X_0' =
\emptyset$. (i) est évident.

Pour $n\ge 1$, on prend une composante \irr $U_1$ de $X_\eta$.
Soit $X_1$ l'adhérence de $U_1$ dans $X$. C'est une composante
\irr de $X$. Soit $X_2$ la réunion des autres composantes \irrs de
$X$. On munit $X_1$ et $X_2$ des structures de \sch réduit
induites. $(X_2)_\eta$ a $n-1$ composantes \irrs, qui sont \geot
\irrs. On a un morphisme fini \surj $ X_1 \coprod X_2 \to X$. On
applique (i) à $X_1$ et obtient $S_1 \to S$ et $(X_1)_0'\to
(X_1)_{S_1}$. Il suffit alors d'appliquer l'\hyp de \rec à
$(X_2)_{S_1}$.

(ii) On fait une \rec sur $n$. Quand $n=0$, (ii) dégénère en (i).
Supposons donnés $S_n \to S$ et $r\sm^\np : X^\np\sm \to X_{S_n}$
\ver ant (ii). On applique (i) au \sch
$(\cosq_n(X^\np\sm/X_{S_n}))_{n+1}$ sur $S_n$ (avec $U$ remplacé par
son image inverse) et obtient un morphisme fini de traits $S'\to
S_n$ et un morphisme $\beta : N \to
(\cosq_n((X^\np\sm)_{S'}/X_{S'}))_{n+1}$ propre avec $\beta_\eta$
surjectif \ver ant \ref{cond.superdJ}. Alors le $X_{S'}$-\sch
simplicial $(n+1)$-tronqué \sscin $r\sm : X'\sm \to X_{S'}$ associé
au triplet $((X^\np\sm)_{S'}, N, \beta)$ \ver e les conditions de
(ii) pour $n+1$.
\end{proof}

\begin{proof}[Démonstration de \ref{theo.RP}]
La démonstration est parallèle à celle de \eqref{eq.Rf*1}. 

On fait une \rec sur $d=\dim X_\eta$. Le cas $d<0$ est trivial.

Soit $d\ge 0$. Il faut montrer que pour $\cG\in \Mod_c(X_\eta,
\Qlb)$ \ent (\resp \entinv), $R\Psi_X \cG$ est \ent (\resp
\eentinv[1,0]).

\Ops $X_\eta$ réduit. \Ops $X$ affine, donc \sep. Choisissons $j:
U\hookrightarrow X_\eta$ ouvert \reg \tq $\cG|U$ soit lisse et que
son \compl $Z= X_\eta -U$ soit de dimension $<d$. Soient $\lbar{Z}$
l'adhérence de $Z$ dans $X$, $i: \lbar{Z} \hra X$. On a le \trdist
\[R\Psi_X i_{\eta*} Ri_\eta^! \cG \to R\Psi_X \cG \to R\Psi_X Rj_* j^* \cG \to.\]
Comme $R\Psi_X i_{\eta*} Ri_\eta^! \simeq i_{s*} R\Psi_{\lbar{Z}}
Ri_\eta^! \cG$ est \ent (\resp \eentinv[1,0]) en vertu de l'\hyp de
\rec, il suffit de voir que $R\Psi_X Rj_* j^* \cG$ est \ent (\resp
\eentinv[1,0]).

Soit $\fH = j^* \cG$. $\fH \simeq (\fH_\cO \ot_\cO E)\ot_E \Qlb$
avec $\fH_\cO$ lisse. Soit $p: U'\to U$ un \rev \etale surjectif qui
trivialise $\fH_\cO\ot_\cO (\cO/\fm)$, où $\fm$ est l'idéal maximal
de $\cO$. $\fH$ est facteur direct de $p_* p^* \fH$, de sorte qu'il
suffit de voir l'intégralité (\resp la $I$-\inte inverse) de
$R\Psi_X R(jp)_* p^* \fH$. On factorise le composé $U'\sto{jp}
X_\eta \hra X$ en $U' \stackrel{j'}{\hra} X'\sto{g} X$ où $j'$ est
une immersion ouverte et $g$ propre.
\[R\Psi_X R(jp)_* p^* \fH \simeq Rg_{s*} R\Psi_{X'} Rj'_{\eta*} p^* \fH.\]

Donc il suffit de \ver er que pour $X$ un \sch \sep \tf sur $S$, $j:
U\hra X_\eta$ un ouvert et $\cG \in \Mod_c(U,\Qlb)$ \ent (\resp
\entinv), $\cG\simeq (\cG_\cO \ot_\cO E)\ot_E \Qlb$ avec
$\cG_\cO\ot_\cO (\cO/\fm)$ constant, on a $R\Psi_X Rj_* \cG$ \ent
(\resp \eentinv[1,0]).

Soit $i\ge 0$. On applique \ref{lemm.superdJ} (ii) à $j$ et
$n=i+1$. On obtient un morphisme fini de traits $f: S'\to S$ et un
carré cartésien de \schs simpliciaux $(i+1)$-tronqués \sscins
\[\xymatrix{U'\sm \ar@{^{(}->}[r] \ar[d]^{s\sm} &X'\sm\ar[d]^{r\sm}\\
U_{S'} \ar@{^{(}->}[r] &X_{S'}}\] où $s\sm$ est un \hyperr propre
$(i+1)$-tronqué et $r_m$ \ver e \ref{cond.superdJ}, $0\le m\le
i+1$. On note des changements de base de $f$ encore par $f$.
\[f^* R\Psi_{X/S} Rj_* \cG \simeq R\Psi_{X_{S'}/S'} f^* Rj_* \cG
\simeq R\Psi_{X_{S'}/S'} Rj_{S'*}\cG_{S'},
\]
donc
\begin{multline*}
\tau_{\le i} f^* R\Psi_{X/S} Rj_* \cG \simeq \tau_{\le i}
R\Psi_{X_{S'}/S'} Rj_{S'*}  R{s\sm}_* s\sm^* \cG_{S'} \\
= \tau_{\le i} R\Psi_{X_{S'}/S'} R{r\sm}_{\eta*} R{j'\sm}_* s\sm^*
\cG_{S'} \simeq\tau_{\le i} R{r\sm}_{s*} R\Psi R{j'\sm}_* s\sm^*
\cG_{S'}.
\end{multline*}
Il suffit de voir que $R\Psi_{X'_m/S'} Rj'_{m*} s_m^* \cG_{S'}$ est
\ent (\resp \eentinv[1,0]), $0\le m\le i$. Il suffit alors
d'appliquer \ref{prop.RPss} (ce qui est licite, car les $s_m^*
\cG_{S'}$ sont modérés).
\end{proof}

\begin{vari}\label{ss.vari.trait}
Soient $S$, $\eta$, $s$, $l$ comme dans \ref{vari.Xse}. Soient $X$
un \sch \tf sur~$S$, $\cF \in \Mod_c(X,\Qlb)$. Pour $r \in \Q$,
$\cF$ est dit $r$-\emph{entier} (\resp $r$-\emph{\entinv{}}) si
$\cF_\eta$ et $\cF_s$ le sont. De \me pour les complexes. On prend
$\delta(X) = \max\{\dim X_\eta+1, \dim X_s\}$, $D_X = R\Hom(-,Ra_X^!
\Qlb(1)[2])$, où $a_X : X\to S$. On a les analogues de
\ref{theo.Rf!b} à \ref{theo.RHom}, \ref{prop.Rf*qf},
\ref{prop.i^!}, en remplaçant $\dim$ par $\delta$, $d_r$ par
\[\max_{y\in |Y_\eta| \cup |Y_s|} \dim f^{-1}(y).\]
On a aussi un
analogue de \ref{prop.Rj*dcn} en ajoutant l'hypothèse que $(X,D)$
est un couple semi-stable sur~$S$.

En effet, \ref{theo.Rf!b} pour $S$ découle trivialement de
\ref{theo.Rf!b} pour $s$ et pour $\eta$. Pour l'analogue de
\ref{prop.Rj*dcn}, soit $\cF\in \Mod_c(U,\Qlb)$ lisse, entier, \mr
sur~$X$. D'après \ref{prop.Rj*dcn} pour $\eta$ et \ref{prop.RPss},
$Ru_* \cF$ et $R\Psi_X Ru_* \cF$ sont $I$-entiers, où $u : U \hra
X_\eta$. Soit $I=\Ker(\Gal(\etab/\eta) \to \Gal(\bar{s}/s))$ le
groupe d'inertie. C'est une extension de $\hat{\Z}_{p'}(1)$ par un
pro-$p$-groupe $P$. La suite spectrale de Hochschild-Serre donne
\[E_2^{p,q} = \fH^p(I,R^q \Psi_X Ru_*\cF) \Rightarrow i^*R^{p+q}j_* \cF,\]
où $i : X_s \to X$. Soit $R^q_t = (R^q\Psi_X Ru_*\cF)^P$. Si
$\sigma$ est un générateur de $\hat{\Z}_{p'}(1)$, on a
\[E_2^{0,q}=\Ker(\sigma-1,R^q_t), \quad
E_2^{1,q}=\mathrm{Coker}(\sigma-1,R^q_t) (-1),
\]
et $E_2^{p,q}=0$
pour $p\neq 0,1$. Donc $Rj_* \cF$ est $I$-\ent.

Les résultats \eqref{eq.Rf*1} et \eqref{eq.Rf*3} pour $S$
découlent de ces résultats pour $s$ et pour $\eta$ et de
\ref{theo.RP}, en imitant les arguments dans \cite[Th.\ finitude,
3.11, 3.12]{SGA4d} comme suit. Le cas spécial $f = j_Y : Y_\eta
\hra Y$ résulte de \ref{theo.RP} et de la suite spectrale de
Hochschild-Serre
\[E_2^{p,q} = \fH^p(I,R^q \Psi_Y - ) \Rightarrow i_Y^*R^{p+q}j_{Y*} -,
\]
où $i_Y : Y_s \to Y$. Traitons le cas général. Soit $L\in
\Dbc(X,\Qlb)$ \ent (\resp $I$-\entinv). Soient $i_X:X_s \to X$,
$j_X:X_\eta \hra X$. On a les triangles distingués
\begin{gather*}
i_{X*}Ri_X^! L \to L \to Rj_{X*}j_X^* L \to,\\
R(fi_X)_* Ri_X^! L \to Rf_* L \to R(fj_X)_* j_X^* L \to.
\end{gather*}
D'après le cas spécial, $Rj_{X*} j_X^* L$ est \ent (\resp
$I$-\entinv), donc $Ri_X^! L$ l'est aussi. Comme $fi_X = i_Y f_s$,
$fj_X = j_Y f_\eta$, on conclut en appliquant \eqref{eq.Rf*1} (\resp
\eqref{eq.Rf*3}) pour $s$ et pour $\eta$ et le cas spécial.

Une fois \eqref{eq.Rf*1} et \eqref{eq.Rf*3} établis pour $S$, on
peut refaire \ref{remq.Ri!} à \ref{prop.i^!} et \ref{ss.demo.RHom},
donnant la démonstration de \ref{theo.Rf*b} à \ref{theo.RHom},
sauf \eqref{eq.Rf!1} et \eqref{eq.Rf!3}. Les résultats
\eqref{eq.Rf!1} et \eqref{eq.Rf!3} pour~$S$ découlent de leurs
analogues pour $s$ et pour $\eta$, en appliquant \ref{remq.Ri!} pour
$S$ et \eqref{eq.demo.Rf!} à $U=Y_\eta$, $W=Y_s$. Notons que
\eqref{eq.D2} et \eqref{eq.D4} pour $S$ peuvent aussi se déduire de
leurs analogues pour $s$ et pour $\eta$, en appliquant
\eqref{eq.demo.D} à $U=X_\eta$, $V=X_s$.
\end{vari}

\begin{vari}\label{ss.vari.Weil}
On peut remplacer les faisceaux usuels partout par des \faisx de
Weil \cite[1.1.10]{WeilII}. Tous les résultats et les variantes qui
précèdent restent valables.
\end{vari}

\newcommand{\variantefausse}{
\begin{vari}
Soient $l$ un nombre premier, $S=\Spec \Z[1/l]$. Soient $X$ un \sch
\tf sur $S$, $\cF \in \Mod_c(X,\Qlb)$. Pour $r \in \Q$, $\cF$ est
dit $r$-\emph{entier} (\resp $r$-\emph{\entinv{}}) si pour tout
point fermé $s$ de $S$, $\cF_s$ l'est. De \me pour les complexes.
On prend $D_X = R\cHom(-,Ra_X^! \Qlb(1)[2])$, où $a_X : X\to S$. On
a les analogues de \ref{theo.Rf!b} à \ref{theo.RHom},
\ref{prop.Rj*dcn}, \ref{prop.Rf*qf}, \ref{prop.i^!}.

En effet, \ref{theo.Rf!b} pour $S$ découle trivialement de
\ref{theo.Rf!b} pour un corps fini. Les autres résultats admettent
les démonstrations similaires à celles données aux \S\S\ 3 et 4.
Dans la démonstration de \ref{lemm.dJ} (i), on remplace
\cite[4.1]{deJong} par \cite[8.2]{deJong}. On utilise aussi les
faits suivants : un \sch \irr \tf sur $S$ est biéquidimensionnel
\cite[10.6.1]{EGAIV} (ce fait n'est utilisé que dans la
démonstration de \ref{theo.Rf!h} pour $S$) ; un \sch $X$ \tf sur $S$
et un ouvert dense $U$ de $X$ ont la \me dimension \ibid[10.6.2].
\end{vari}
}

\begin{vari}\label{ss.vari.A}
Soit $A$ un sous-anneau intégralement fermé  de $\Qlb$. On pose
\[A^{-1} = \ensdr{\al \in \Qlb^\times}{\al^{-1}\in A}.\]
On fixe un plongement $\iota : \Qb \to \Qlb$. Avec les notation du
\S~\ref{S1}, pour $r\in \Q$, un $\Qlb$-\faisc $\cF$ sur un \sch $X$
\tf sur $\eta$ est dit $r$-$A$-\emph{\ent{}} (\resp
$r$-$A$-\emph{\entinv{}}) si pour tout $x\in |X|$, les valeurs
propres de $\Phi_x$ sont dans $\iota(q^{r}) A$ (\resp $\iota(q^{r})
A^{-1}$), où $q=\sharp \kappa(x)$. Cette définition ne dépend pas
des choix de $\Phi_x$ et de $\iota$. On dit que $\cF$ est
$A$-\emph{\ent{}} (\resp $A$-\emph{\entinv{}}) s'il est $0$-$A$-\ent
(\resp $0$-$A$-\entinv). On définit aussi la $A$-intégralité des
complexes. Tous les résultats et les variantes restent valables
pour ces notions.

\vcb{v1m}{Si $A$ est de plus \cic \cite[V, \S~1, \no 4,
déf.~5]{BourAC} (en particulier si $A$ est la fermeture intégrale
d'un sous-anneau noethérien intégralement clos de $\Qlb$
\ibid[exerc. 14]), alors d'après un lemme de Fatou
\cite[8.3]{IllMisc}, $\cF$ est $r$-$A$-\ent \ssi pour tout $x\in
|X|$ et tout entier $n\ge 1$, $\Tr(\Phi_x^n, \cF_{\xb})$ appartient
à $\iota(q^{nr}) A$. Ce critère n'a pas d'analogue pour les
faisceaux entiers inverses.}

Si on prend pour $A$ la fermeture intégrale de $\Z$ dans $\Qlb$,
on retrouve la notion d'intégralité dans ce qui précède.

Soit $T$ un ensemble de nombres premiers. Si on prend pour $A$ la
fermeture intégrale de $\Z[1/t]_{t\in T}$ dans $\Qlb$, on retrouve
la notion de $T$-intégralité dans \cite[XXI~5]{SGA7} et \cite{DE}.
\end{vari}

\newcommand{\ancienparagraphesix}{
\section{Appendice : Raffinements}
Le but de cet appendice est d'éliminer l'\hyp de séparation dans
\eqref{eq.Rf*2}. On a besoin d'une analyse plus raffinée de la
situation.

On reprend les notations du \S~1. Soient $X$ un \sch \tf sur $\eta$,
$\cF$ un $\Qlb$-\faisc sur $X$.

\begin{defn}\label{defn.fe}
Soient $e\in \N$, $r\in \Q$. On dit qu'un $\Qlb$-\faisc $\cF$ est
$r$-\emph{entier} (\resp $r$-\emph{\entinv{}}) en dimension $\ge e$
s'il existe une partie fermée $Y\hra X$ de dimension $<e$ \tq
$\cF|(X-Y)$ soit \ent[r] (\resp \entinv[r]). On dit que $\cF$ est
\emph{entier} (\resp \emph{entier inverse}) en dimension $\ge e$
s'il est $0$-entier (\resp $0$-\entinv) en dimension $\ge e$.
\end{defn}

Si $e = 0$, \og \ent[r] (\resp \entinv[r]) en dimension $\ge e$ \fg
équivaut à \og \ent[r] (\resp \entinv[r]) \fg.

Soit $e\in \N$. Les $\Qlb$-\faisx \ents (\resp $r$-entiers, \resp
\entinvs, \resp $r$-\entinvs) en dimension $\ge e$ sur $X$ forment
une sous-\cat épaisse de $\Mod_c(X,\Qlb)$. On la désigne par
$\dModclent{X}{e}$ (\resp $\dModclent[r]{X}{e}$, \resp
$\dModclentinv{X}{e}$, \resp $\dModclentinv[r]{X}{e}$).

\begin{defn}\label{defn.ce}
Soit $\eps: \Z \to \Q$ une fonction. On dit qu'un objet $K\in
\Dbc(X,\Qlb)$ est \emph{entier} (\resp $\eps$-\emph{entier}, \resp
\emph{entier inverse}, \resp $\eps$-\emph{\entinv{}}) en dimension
$\ge e$ si pour tout $i\in \Z$, $\fH^i(K)$ est entier (\resp
$\eps(i)$-entier, \resp \entinv, \resp $\eps(i)$-\entinv) en
dimension $\ge e$.
\end{defn}

On désigne la sous-\cat pleine de $\Dbc(X,\Qlb)$ formée des objets
\ents (\resp \eents, \resp \entinvs, \resp \eentinvs) en dimension
$\ge e$ par
\[\text{$\dDbclent{X}e$ (\resp $\dDbcleent{X}e$, \resp $\dDbclentinv{X}e$,
\resp $\dDbcleentinv{X}e$).}\]

On va énoncer dans ce cadre des raffinements des théorèmes
\ref{theo.Rf!b} à \ref{theo.D}.

\begin{prop}\label{theo.Rf!be} Soient $f: X\to Y$ un morphisme \sep de \schs \tf
sur~$\eta$, $d_r = \max_{y\in |Y|} \dim f^{-1}(y)$. Alors $Rf_!$
induit
\begin{align}
\dDbcent{X}e &\to \dDbcent{Y}e,
\\
\dDbceent[1,0]{X}e &\to \dDbceent[1,-d_r]{Y}e,
\\
\dDbceentinv[1,0]{X}e &\to
\dDbceentinv[1,0]{Y}e,
\\
\dDbcentinv{X}e &\to \dDbcentinv[d_r]{Y}e.
\end{align}
\end{prop}

\begin{prop}\label{theo.Rf*be}
Soient $f: X\to Y$ un morphisme de \schs \tf sur $\eta$, $d_X=\dim
X$. Alors $Rf_*$ induit
\begin{align}
\dDbcent{X}e &\to \dDbcent{Y}e,\label{eq.Rf*1e}\\
\dDbceent[1,0]{X}e &\to \dDbceent[1,-d_X+e]{Y}e,\label{eq.Rf*2e}\\
\dDbceentinv[1,0]{X}e &\to \dDbceentinv[1,0]{Y}e,\label{eq.Rf*3e}\\
\dDbcentinv{X}e &\to \dDbcentinv[(d_X-e)]{Y}e.\label{eq.Rf*4e}
\end{align}
\end{prop}

\Enpart \ref{theo.Rf*b} est vrai sans hypothèse de séparation.

\begin{prop}\label{theo.Rf!he}
Soient $f: X\to Y$ un morphisme \sep de \schs \tf sur $\eta$,
$d_Y=\dim Y$, $d_r = \max_{y\in |Y|} \dim f^{-1}(y)$
. Alors $Rf^!$ induit
\begin{align}
\dDbcent{Y}{e-d_r} &\to \dDbcent[-d_r]{X}e,\label{eq.Rf!1e}\\
\dDbceent[1,0]{Y}{e-d_r} &\to \dDbceent[1,-d_Y+e]{X}e,\label{eq.Rf!2e}\\
\dDbceentinv[1,0]{Y}{e-d_r} &\to \dDbceentinv[1,d_r]{X}e,\label{eq.Rf!3e}\\
\dDbcentinv{Y}{e-d_r} &\to
\dDbcentinv[(d_Y-e)]{X}e.\label{eq.Rf!4e}
\end{align}
\end{prop}

\begin{prop}\label{theo.De}
Soient $X$ un \sch \tf sur $\eta$, $d_X=\dim X$. Alors $D_X$ induit
\begin{align}
\dDbcent{X}e\op &\to \dDbcentinv[-e]{X}e,\label{eq.D1e}\\
\dDbceent[1,0]{X}e\op &\to \dDbceentinv[1,d_X]{X}e,\label{eq.D2e}\\
\dDbceentinv[1,0]{X}e\op &\to \dDbceent[1,e]{X}e,\label{eq.D3e}\\
\dDbcentinv{X}e\op &\to \dDbcent[-d_X]{X}e,\label{eq.D4e}
\end{align}
De plus, pour $K\in\dModclentinv{X}e$, $\fH^a(DK)$ est \ent[(a+1)]
en \dim $\ge e$, $-d_X\le a \le -e-1$.
\end{prop}


\begin{proof}[\Demo de \ref{theo.Rf!be} à \ref{theo.De}
] Résultent facilement des résultats analogues du \S~2 :
\ref{theo.Rf!be}, \eqref{eq.Rf*1e}, \eqref{eq.Rf*3e},
\eqref{eq.Rf!1e}, \eqref{eq.Rf!3e}, \eqref{eq.D2e},
\eqref{eq.D4e}
. Le reste, sauf \eqref{eq.Rf*2e} pour $f$ non séparé, admet les
mêmes démonstrations que celles données au \S~4.

Il reste à montrer \eqref{eq.Rf*2e} pour $f$ non séparé. Soit
$K\in \dModcent{X}e$. On va montrer que pour tout sous-\sch fermé
$i : W\hra X$, on a $Rf_* i_* Ri^! K \in \dDbcent[(I-d_X+e)]{X}e$.
On fait une \rec sur $d_W =\dim W$. Le cas $d_W < e$ est trivial.
Pour $d_W\ge e$, on a $K\in \dModcent{X}{d_W}$, donc $Ri^! K \in
\dDbcent[(I-d_X+d_W)]{W}{d_W}$ en vertu de \eqref{eq.Rf!2e}. Par
définition, il existe un ouvert $j: U\hra W$ de \compl $i_Z : Z\hra
W$ de dimension $<d_W$ \tq $(Ri^!K)|U$ soit $(I-d_X + d_W)$-\ent.
Quitte à rétrécir $U$, \ops $U$ séparé, voire affine. On a le
\trdist
\[Rf_*i_*i_{Z*}Ri_Z^! Ri^!K \to Rf_*i_*Ri^!K \to Rf_*i_* Rj_*j^*Ri^!K \to.\]
D'après \eqref{eq.Rf*2e} (cas séparé), on a
\[Rf_*i_* Rj_*j^*Ri^!K =
R(fij)_*((Ri^!K)|U) \in \dDbcent[(I-d_X+e)]{Y}e.\] \Dpr l'\hyp de
\rec, on a
\[Rf_*i_*i_{Z*}Ri_Z^! Ri^!K= Rf_* (ii_Z)_* R(ii_Z)^! K \in
\dDbcent[(I-d_X+e)]{Y}e.\] D'où le résultat.
\end{proof}

On laisse au lecteur le soin de donner des raffinements des
 autres résultats.
}

\section{Appendice : Intégralité sur les champs
algébriques}\label{S6}

\newcommand{\chU}{{\mathcal{U}}}
\newcommand{\cX}{{\mathcal{X}}}
\newcommand{\cY}{{\mathcal{Y}}}
\newcommand{\cZ}{{\mathcal{Z}}}

\newcommand{\ment}{{\mathrm{ent}}}
\newcommand{\mentinv}{{\mathrm{ent}^{-1}}}
\newcommand{\mtent}{\textrm{-}\ment}
\newcommand{\mtentinv}{\textrm{-}\mentinv}

\newcommand{\DM}{de Deligne-Mumford\xspace}

Soient $K$ un corps fini de \cara $p$ ou un corps local de \cara
résiduelle~$p$, $\eta=\Spec K$, $l$ un nombre premier $\neq p$.
Pour $\cX$ un $\eta$-champ algébrique \cite[4.1]{Laumon-MB} \tf, on
note $\Mod_c(\cX,\Qlb)$ la catégorie des $\Qlb$-faisceaux
constructibles sur le site lisse-étale de $\cX$ \ibid[12.1 (i)]. On
dispose, par \cite{Laszlo-Olsson2}, d'une catégorie triangulée
$D_c(\cX,\Qlb)$ munie d'une $t$-structure de c\oe ur
$\Mod_c(\cX,\Qlb)$. On écrira $\Mod_c(\cX)$ pour $\Mod_c(\cX,\Qlb)$
et $D_c(\cX)$ pour $D_c(\cX,\Qlb)$. On dispose d'un formalisme de
six opérations : pour $f: \cX\to \cY$ un morphisme de $\eta$-champs
algébriques \tf,
\begin{align*}
D_\cX : D_c(\cX)\op &\to D_c(\cX),\\
-\otimes- : D_c^-(\cX) \times D_c^-(\cX) &\to D_c^-(\cX),\\
R\cHom_\cX(-,-) : D_c^-(\cX)\op \times D_c^+(\cX) &\to D_c^+(\cX),\\
Rf_* : D_c^+(\cX) &\to D_c^+(\cY),\\
Rf_! : D_c^-(\cX) &\to D_c^-(\cY),\\
f^*, Rf^! : D_c(\cY) &\to D_c(\cX).
\end{align*}
Si $\cX$ est un $\eta$-champ \DM \tf, $\Mod_c(\cX)$ s'identifie à
la \cat des $\Qlb$-faisceaux constructibles sur le site
\emph{étale} de $\cX$ \cite[12.1 (ii)]{Laumon-MB}.

\begin{defn}
Soit $\eps : \Z\to \Q$ une fonction. On dit que $L\in D_c(\cX)$ est
$\eps$-\emph{entier} (\resp $\eps$-\emph{entier inverse}) si pour
tout point $i: x\to \cX$ avec $\kappa(x)$ une extension finie de~$K$
et tout $a\in \Z$, $\fH^a(i^* L) \in \Mod_c(x,\Qlb)$ est
$\eps(a)$-entier (\resp $\eps(a)$-\entinv) (au sens de
\ref{defn.f}).
\end{defn}

Si $\cX$ est un schéma \tf sur $\eta$ et $L\in \Dbc$, alors cette
définition coïncide avec \ref{defn.c}.

Soit $f: \cX\to \cY$ un morphisme de $\eta$-champs algébriques de
type fini. Si $M\in D_c(\cY)$ est \eent (\resp \eentinv), il en est
de même de $f^*M\in D_c(\cX)$. La réciproque est vraie lorsque $f$
est surjectif. Cela donne le critère suivant : $L\in \Dbc(\cX)$ est
\eent (\resp \eentinv) si et seulement si pour tout (ou pour un)
morphisme surjectif $g: X\to \cX$ avec $X$ un schéma de type fini
sur~$\eta$, $g^*L$ est \eent (\resp \eentinv) (au sens de
\ref{defn.c}). Pour $r\in \Q$, $-\otimes-$ induit
\begin{align*}
D_c^-(\cX)_{rI\mtent} \times D_c^-(\cX)_{rI\mtent} \to
D_c^-(\cX)_{rI\mtent}, \\
D_c^-(\cX)_{rI\mtentinv} \times D_c^-(\cX)_{rI\mtentinv} \to
D_c^-(\cX)_{rI\mtentinv}.
\end{align*}
Pour $\cF,\cG\in \Mod_c(\cX)$ avec $\cF$ lisse, si $\cF$ est \entinv
(\resp entier) et $\cG$ est entier (\resp \entinv), on a
\[\cHom(\cF,\cG) \in \Mod_c(\cX)_{\ment} \text{ (\resp $\in \Mod_c(\cX)_{\mentinv}$)}.\]

On suppose dorénavant que $\cX$ est non vide. Rappelons qu'une
présentation $P: X\to \cX$ est un morphisme surjectif lisse avec
$X$ un espace algébrique. On pose $c_\cX = \min_P \dim P \in \N$,
où $P$ parcourt les présentations $P : X\to \cX$, $\dim P =
\sup_{x\in X} \dim_x P$ \cite[p.~98]{Laumon-MB}, $d_\cX = \dim \cX
\in \Z$ \ibid[(11.15)]. Par définition, $d_\cX \ge -c_\cX$. On a
$c_{\cX}=0$ si et seulement si $\cX$ est un $\eta$-champs de
Deligne-Mumford.
On pose $c_r = \min \dim P \in \N$, où
le minimum est pris sur tous les systèmes
\[\xymatrix{X\ar[r]^P & \cX\times_\cY Y\ar[r]\ar[d]
& \cX\ar[d]^{f}\\
& Y \ar[r]^Q & \cY}
\]
où $P$ et $Q$ sont des présentations et le carré est
$2$-cartésien, $d_r = \dim f = \max_\xi \dim \cX_\xi \in \Z$, où
$\xi$ parcourt les points de~$\cY$. On a $d_r\ge -c_r$. Rappelons
que $f$ est dit \emph{relativement de Deligne-Mumford} \ibid[7.3.3]
si pour tout schéma affine $Y$ et tout morphisme $Y\to \cY$, le
produit fibré $Y\times_{\cY} \cX$ est un $\eta$-champ de
Deligne-Mumford. On a $c_r=0$ si et seulement si $f$ est un
morphisme relativement de Deligne-Mumford. On a $c_r \le c_{\cX} \le
c_{\cY} + c_r$, $d_r-c_{\cY} \le d_{\cX} \le d_{\cY} + d_r$.

\begin{prop}\label{theo.Dch}
Le foncteur $D_\cX$ induit
\begin{align*}
D_c(\cX)_\mentinv\op &\to D_c(\cX)_{-d_\cX\mtent},\\
D_c(\cX)_{I\mtentinv}\op &\to D_c(\cX)_{(I-c_\cX)\mtent},\\
D_c(\cX)_{I\mtent}\op &\to D_c(\cX)_{(I+d_\cX)\mtentinv},\\
D_c(\cX)_\ment\op &\to D_c(\cX)_{c_\cX\mtentinv}.
\end{align*}
De plus, pour $\cF\in \Mod_c(X)_{\mentinv}$, $\fH^a(D_\cX \cF)$ est
$(a-c_\cX +1)$-entier, $a\le c_\cX-1$.
\end{prop}

\begin{proof}
On prend une présentation $P: X \to \cX$ purement de dimension
$c_\cX$ avec $X$ un \sch \tf sur $\eta$. On a $\dim X \le d_\cX +
c_\cX$. Pour $L\in D_c(\cX)$, $P^* D_\cX L \simeq D_X RP^! L \simeq
(D_X P^* L) (-c_\cX)[-2c_\cX]$. Comme l'amplitude cohomologique de
$D_X$ est bornée, il suffit donc d'appliquer \ref{theo.D}.
\end{proof}

\begin{prop}
Le foncteur $Rf^!$ induit
\begin{align*}
D_c(\cY)_\ment &\to D_c(\cX)_{-d_r\mtent},\\
D_c(\cY)_{I\mtent} &\to D_c(\cX)_{(I-d_\cY-c_\cY-c_r)\mtent},\\
D_c(\cY)_{I\mtentinv} &\to D_c(\cX)_{(I+d_r)\mtentinv},\\
D_c(\cY)_{\mentinv} &\to D_c(\cX)_{(d_\cY+c_\cY+c_r)\mtentinv}.
\end{align*}
\end{prop}

\begin{proof}
Formons le diagramme à carré $2$-cartésien
\[\xymatrix{X\ar[r]^P & \cX\times_\cY Y\ar[r]^{Q'}\ar[d]^{f_Y}
& \cX\ar[d]^f\\
& Y \ar[r]^Q & \cY}
\]
où $Q$ est une présentation purement de dimension $c_\cY$, $P$ est
une présentation purement de dimension $c_r$, $Y$ est un \sch
quasi-compact, $X$ est un \sch affine (donc séparé sur~$Y$). On a
$\dim Y \le d_\cY + c_\cY$, $\dim (f_Y \circ P) \le d_r + c_r$. Pour
$L\in D_c(\cY)$,
\[(Q' \circ P)^*Rf^!L \simeq P^* Rf_Y^! Q^* L
\simeq R(f_Y\circ P)^! Q^*L (-c_r)[-2c_r].
\]
Il suffit alors
d'appliquer \ref{theo.Rf!h}.
\end{proof}

\begin{prop}\label{prop.Rf*Rf!ch}
Le foncteur $Rf_*$ induit
\begin{align}
D_c^+(\cX)_\ment &\to D^+_c(\cY)_\ment,\label{eq.Rf*1ch}\\
D_c^+(\cX)_{I\mtentinv} &\to
D^+_c(\cY)_{I\mtentinv}\label{eq.Rf*3ch},
\end{align}
et $Rf_!$ induit
\begin{align}
D_c^-(\cX)_{I\mtent} &\to D_c^-(\cY)_{(I-d_r)\mtent},\label{eq.Rf!b2ch}\\
D_c^-(\cX)_{\mentinv} &\to
D^-_c(\cY)_{d_r\mtentinv}\label{eq.Rf!b4ch}.
\end{align}
Si $f$ est relativement \DM, $Rf_*$ induit
\begin{align}
D_c^+(\cX)_{I\mtent} &\to D_c^+(\cY)_{(I-d_\cX-c_\cY)\mtent},\label{eq.Rf*2ch}\\
D_c^+(\cX)_{\mentinv} &\to
D_c^+(\cY)_{(d_\cX+c_\cY)\mtentinv},\label{eq.Rf*4ch}
\end{align}
et $Rf_!$ induit
\begin{align}
D_c^-(\cX)_{\ment} &\to D_c^-(\cY)_{\ment},\label{eq.Rf!b1ch}\\
D_c^-(\cX)_{I\mtentinv} &\to
D_c^-(\cY)_{I\mtentinv}.\label{eq.Rf!b3ch}
\end{align}
\end{prop}

\begin{proof}
Soit $Y \to \cY$ une présentation avec $Y$ est un schéma \tf
sur~$\eta$. Pour \eqref{eq.Rf*1ch} et \eqref{eq.Rf*3ch}, quitte à
remplacer $\cY$ par $Y$, \ops que $\cY$ est un schéma. On prend un
hyperrecouvrement lisse $P\sm : X\sm \to \cX$ où les $X_n$ sont des
\schs affines (donc séparés sur $\cY$). Pour $L\in D^+_c(\cX)$,
$Rf_* \simto Rf_* R{P\sm}_* P\sm^* L$. Il suffit alors d'appliquer
\eqref{eq.Rf*1} et \eqref{eq.Rf*3}.

Pour les résultats concernant $Rf_!$, \ops que $\cY$ est le spectre
d'un corps. On a $d_\cX = d_r$. Alors \eqref{eq.Rf!b2ch} et
\eqref{eq.Rf!b4ch} découlent du dernier alinéa et de
\ref{theo.Dch} : $Rf_! = D_\cY Rf_* D_\cX$ induit
\begin{gather*}
D_c^-(\cX)_{I\mtent} \sto{D_\cX} D_c^+(\cX)_{(I+d_\cX)\mtentinv}\op
\sto{Rf_*} D_c^+(\cY)_{(I+d_\cX)\mtentinv}\op \sto{D_\cY}
D_c^-(\cY)_{(I-d_\cX)\mtent},\\
D_c^-(\cX)_{\mentinv} \sto{D_\cX} D_c^+(\cX)_{-d_\cX\mtent}\op
\sto{Rf_*} D_c^+(\cY)_{-d_\cX\mtent}\op \sto{D_\cY}
D_c^-(\cY)_{d_\cX\mtent}.
\end{gather*}

Pour \eqref{eq.Rf!b1ch} et \eqref{eq.Rf!b3ch}, on est donc ramené
au cas où $\cX$ est un $\eta$-champ de \DM. On fait une récurrence
sur $d_\cX$. Il existe une immersion ouverte dominante $j : \chU
\hra \cX$ et un morphisme fini étale $\pi : U \to \chU$, où $U$
est un \sch affine \cite[6.1.1]{Laumon-MB}. Soient $\cZ$ le fermé
complémentaire de $\chU$ dans~$\cX$, $i: \cZ \to \cX$. Pour $L\in
D_c^-(\cX)$, le \trdist
\[j_! j^* L \to L \to i_* i^* L \to\]
induit le triangle distingué
\[R(fj)_! j^* L \to Rf_! L \to R(fi)_! i^* L \to.\]
Comme $j^* L$ est facteur direct de $\pi_* \pi^* j^*L$, il suffit
d'appliquer \eqref{eq.Rf!b1} et \eqref{eq.Rf!b3} à $R(fj\pi)_!
(j\pi)^*L$ et l'hypothèse de récurrence à $R(fi)_! j^*L$.

Enfin, \eqref{eq.Rf*2ch} et \eqref{eq.Rf*4ch} résultent du dernier
alinéa et de \ref{theo.Dch} : $Rf_* \simeq D_\cY Rf_! D_\cX$ induit
\begin{gather*}
D_c^+(\cX)_{I\mtent} \sto{D_\cX} D_c^-(\cX)_{(I+d_\cX)\mtentinv}\op
\sto{Rf_!} D_c^-(\cY)_{(I+d_\cX)\mtentinv}\op \sto{D_\cY}
D_c^+(\cY)_{(I-d_\cX-c_\cY)\mtent},\\
D_c^+(\cX)_{\mentinv} \sto{D_\cX} D_c^-(\cX)_{-d_\cX\mtent}\op
\sto{Rf_!} D_c^-(\cY)_{-d_\cX\mtent}\op \sto{D_\cY}
D_c^+(\cY)_{(d_\cX+c_\cY)\mtent}.
\end{gather*}
\end{proof}

\begin{coro}\label{coro.Rf*2}
L'assertion \eqref{eq.Rf*2} est vraie sans hypothèse de
séparation.
\end{coro}

\begin{proof}
C'est un cas particulier de \eqref{eq.Rf*2ch}.
\end{proof}

\begin{remqe}
(i) Le premier alinéa de la démonstration de \ref{prop.Rf*Rf!ch}
montre que $Rf_*$ envoie $\Mod_c(\cX)_\mentinv$ dans
$D_c^+(\cY)_{\eps\mtentinv}$, où
\[\eps(a) =
\begin{cases}
a & \text{si $0\le a \le d_\cX + c_\cX + c_\cY$},\\
a-E\left(\frac{a-d_\cX -c_\cY}{c_\cX +1}\right) & \text{si $a\ge
d_\cX+ c_\cX + c_\cY$}.
\end{cases}
\]
Ici $E$ est la fonction partie entière. Lorsque $f$ n'est pas
relativement \DM, ceci légèrement améliore \eqref{eq.Rf*3ch}. On
peut en déduire une légère amélioration de \eqref{eq.Rf!b2ch}.

(ii) Si $f$ est un morphisme séparé, représentable et quasi-fini \cite[3.10.1]{Laumon-MB}
avec $d_\cX + c_\cY \ge 1$, on a un analogue de
\ref{prop.Rf*qf} qui améliore \eqref{eq.Rf*2ch} : $Rf_*$ envoie
$D_c^+(\cX)_{I\mtent}$ dans $D_c^+(\cY)_{(I+1-d_\cX-c_\cY)\mtent}$.
\end{remqe}

\begin{prop} Le foncteur $R\cHom(-,-)$ induit
\begin{align*}
D_c^-(\cX)_\mentinv\op \times D_c^+(\cX)_\ment &\to
D_c^+(\cX)_\ment,\\
D_c^-(\cX)_{I\mtentinv}\op \times D_c^+(\cX)_{I\mtent} &\to
D_c^+(\cX)_{(I-d_\cX-c_\cX)\mtent},\\
D_c^-(\cX)_{I\mtent}\op \times D_c^+(\cX)_{I\mtentinv} &\to
D_c^+(\cX)_{I\mtentinv},\\
D_c^-(\cX)_\ment\op \times D_c^+(\cX)_\mentinv &\to
D_c^+(\cX)_{(d_\cX+c_\cX)\mtentinv}.
\end{align*}
\end{prop}

\begin{proof}
On procède comme en \ref{ss.demo.RHom}.
\end{proof}

On peut aussi considérer l'intégralité sur les champs
algébriques
sur un trait excellent de corps résiduel fini, ce qui généralise
\ref{ss.vari.trait}. Les résultats sont similaires à ceux exposés
dans ce \S, avec des modifications appropriées des estimations de
dimension.

Les variantes \ref{ss.vari.Weil} et \ref{ss.vari.A} restent toujours
valables.

\renewcommand{\refname}{Bibliographie}
\bibliographystyle{smfalpha}
\bibliography{inte2c}

\providecommand{\SortNoop}[1]{}
\providecommand{\bysame}{\leavevmode ---\ }
\providecommand{\og}{``}
\providecommand{\fg}{''}
\providecommand{\smfandname}{et}
\providecommand{\smfedsname}{\'eds.}
\providecommand{\smfedname}{\'ed.}
\providecommand{\smfmastersthesisname}{M\'emoire}
\providecommand{\smfphdthesisname}{Th\`ese}
\begin{thebibliography}{{\SortNoop{zz}}{\SortNoop{}}{\sigle{EGAIV}}67}

\bibitem[BBD82]{BBD}
{\scshape A.~A. Be{\u\i}linson, J.~Bernstein {\normalfont \smfandname}
  P.~Deligne} -- {\og Faisceaux pervers\fg}, \emph{Analyse et topologie sur les
  espaces singuliers (I)}, Astérisque, vol. 100, Soc.\ math.\ France, 1982,
  p.~5--171.

\bibitem[BO78]{BO}
{\scshape P.~Berthelot {\normalfont \smfandname} A.~Ogus} -- \emph{\eng{Notes
  on crystalline cohomology}}, \eng{Princeton Univ.\ Press}, 1978.

\bibitem[Bou85]{BourAC}
{\scshape N.~Bourbaki} -- \emph{Algèbre commutative, chapitres {V} à {VII}},
  Éléments de mathématique, Masson, 1985.

\bibitem[DE06]{DE}
{\scshape P.~Deligne {\normalfont \smfandname} H.~Esnault} -- {\og
  \eng{Appendix to ``Deligne's integrality theorem in unequal characteristic
  and rational points over finite fields''}\fg}, \emph{Ann.\ Math.\ (2)}
  \textbf{164} (2006), p.~726--730.

\bibitem[Del80]{WeilII}
{\scshape P.~Deligne} -- {\og La conjecture de {Weil : II}\fg}, \emph{Publ.\
  math.\ IHÉS} \textbf{52} (1980), p.~137--252.

\bibitem[dJ96]{deJong}
{\scshape A.~J. de~Jong} -- {\og \eng{Smoothness, semi-stability and
  alterations}\fg}, \emph{Publ.\ math.\ IHÉS} \textbf{83} (1996), p.~51--93.

\bibitem[Eke90]{Ekedahl}
{\scshape T.~Ekedahl} -- {\og \eng{On the adic formalism}\fg}, The Grothendieck
  Festschrift, vol.\ II, Progr.\ Math., vol.~87, Birkhäuser, 1990, p.~197--218.

\bibitem[Fuj02]{Fuji}
{\scshape K.~Fujiwara} -- {\og \eng{A proof of the absolute purity conjecture
  (after Gabber)}\fg}, \emph{\eng{Algebraic geometry 2000, Azumino}}, Adv.
  Stud. Pure Math., vol.~36, Math.\ Soc.\ Japan, 2002, p.~153--183.

\bibitem[Gro57]{Tohoku}
{\scshape A.~Grothendieck} -- {\og Sur quelques points d'algèbre
  homologique\fg}, \emph{Tôhoku Math.\ J. (2)} \textbf{9} (1957), p.~119--221.

\bibitem[Ill02]{IllPL}
{\scshape L.~Illusie} -- {\og Sur la formule de {Picard-Lefschetz}\fg},
  \emph{\eng{Algebraic geometry 2000, Azumino}}, Adv.\ Stud.\ Pure Math.,
  vol.~36, Math.\ Soc.\ Japan, 2002, p.~249--268.

\bibitem[Ill04]{Illss}
\bysame , {\og \eng{On semistable reduction and the calculation of nearby
  cycles}\fg}, \emph{\eng{Geometric aspects of Dwork theory}}, vol.~II, Walter
  de Gruyter, 2004, p.~785--803.

\bibitem[Ill06]{IllMisc}
\bysame , {\og \eng{Miscellany on traces in $l$-adic cohomology: a survey}\fg},
  \emph{Japanese J.~Math.\ (3)} \textbf{1} (2006), p.~107--136.

\bibitem[LMB00]{Laumon-MB}
{\scshape G.~Laumon {\normalfont \smfandname} L.~Moret-Bailly} -- \emph{Champs
  algébriques}, Springer, 2000.

\bibitem[LO06]{Laszlo-Olsson2}
{\scshape Y.~Laszlo {\normalfont \smfandname} M.~Olsson} -- {\og \eng{The six
  operations for sheaves on Artin stacks II: adic coefficients}\fg},
  \eng{\texttt{arXiv:math/0603680v1}}, 2006.

\bibitem[Org03]{Orgogozo}
{\scshape F.~Orgogozo} -- {\og Altérations et groupe fondamental premier
  à~{$p$}\fg}, \emph{Bull.\ Soc.\ math.\ France} \textbf{131} (2003),
  p.~123--147.

\bibitem[Ser98]{Serre}
{\scshape J.-P. Serre} -- \emph{Représentations linéaires des groupes finis},
  5\ieme \smfedname, Hermann, 1998.

\bibitem[Zhe07]{Zheng}
{\scshape W.~Zheng} -- {\og Sur l'indépendance de {$l$} en cohomologie
  {$l$}-adique sur les corps locaux\fg}, \eng{\texttt{arXiv:0711.3658v1}},
  2007.

\bibitem[{\SortNoop{zz}}{\SortNoop{}}{\sigle{EGAIV}}67]{EGAIV}
{\scshape A.~{\SortNoop{zz}}{\SortNoop{}}{\sigle{EGAIV}}Grothendieck} -- {\og
  Éléments de géométrie algébrique : {IV. É}tude locale des schémas et des
  morphismes de schémas\fg}, \emph{Publ.\ math.\ IHÉS} \textbf{20, 24, 28, 32}
  (1964--1967).

\bibitem[{\SortNoop{zz}}{\SortNoop{}}{\sigle{SGA1}}03]{SGA1}
\emph{Revêtements étales et groupe fondamental} -- Séminaire de géométrie
  algébrique du Bois-Marie 1960--1961, dirigé par A.~Grothendieck, Documents
  math., vol.~3, Soc.\ math.\ France, 2003.

\bibitem[{\SortNoop{zz}}{\SortNoop{}}{\sigle{SGA4}}73]{SGA4}
\emph{Théorie des topos et cohomologie étale des schémas} -- Séminaire de
  géométrie algébrique du Bois-Marie 1963--1964, dirigé par M.~Artin,
  A.~Grothendieck, \mbox{J.-L.}~Verdier, LNM, vol. 269, 270, 305,
  Springer-Verlag, 1972--1973.

\bibitem[{\SortNoop{zz}}{\SortNoop{}}{\sigle{SGA4\textonehalf}}77]{SGA4d}
{\scshape P.~{\SortNoop{zz}}{\SortNoop{}}{\sigle{SGA4\textonehalf}}Deligne} --
  \emph{Cohomologie étale}, LNM, vol. 569, Springer-Verlag, 1977.

\bibitem[{\SortNoop{zz}}{\SortNoop{}}{\sigle{SGA7}}73]{SGA7}
\emph{Groupes de monodromie en géométrie algébrique} -- Séminaire de géométrie
  algébrique du Bois-Marie 1967--1969, I, dirigé par A.~Grothendieck, II, par
  P.~Deligne, N.~Katz, LNM, vol. 288, 340, Springer-Verlag, 1972--1973.

\end{thebibliography}
\end{document}